\documentclass[%
 print,
 superscriptaddress,
 amsmath,amssymb,
aps,
]{revtex4-2}

\usepackage{graphicx}
\usepackage{dcolumn}
\usepackage{bm}

\renewcommand{\b}[1]{\boldsymbol{#1}}
\usepackage{color}
\usepackage{amsmath}
\usepackage{algorithm, algpseudocode}
\usepackage{siunitx}

\begin{document}

\title{Interpretable low-order representation of eigenmode\\ deformation in parameterized dynamical systems}

\author{Nicol\'as Torres-Ulloa}
\affiliation{Department of Mechanical Engineering, University of Chile, Beauchef 851, Santiago, Chile}

\author{Erick Kracht}%
\affiliation{Department of Mechanical Engineering, University of Chile, Beauchef 851, Santiago, Chile}

\author{Urban Fasel}%
\affiliation{Department of Aeronautics, Imperial College London, SW7 2AZ, United Kingdom}

\author{Benjamin Herrmann}%
 \email{benjaminh@uchile.cl}
\affiliation{Department of Mechanical Engineering, University of Chile, Beauchef 851, Santiago, Chile}


\begin{abstract}
    Modal analysis has long been consolidated as a basic tool to interpret dynamics and build low-order models of mechanical, thermal, and fluid systems.
    Eigenmodes arising from the spectral decomposition of the underlying linearized dynamics represent spatial patterns in vibration, temperature, or velocity fields associated with simple time dynamics.
    However, for systems that depend on one or more parameters, eigenmodes obtained for one set of parameter values are not necessarily dynamically relevant in other regions of parameter space.
    In this work, we formulate a method to obtain an optimal orthogonal basis of \emph{eigen-deformation modes} (EDMs) that capture eigenmode variations across a range of parameter values.
    Through numerical examples of common parameterized dynamical systems in engineering, we show that EDMs are useful for parameterized model reduction and to provide physical insight into the effects of parameter changes on the underlying dynamics.
\end{abstract}

\maketitle


\section{\label{sec:intro}Introduction}






Understanding and modeling the complex spatiotemporal dynamics of structural vibrations, heat transfer, and fluid flow phenomena is critical for various applications in mechanical, civil, and aerospace engineering. The behavior of these complex systems often depends on quantities such as material or media properties, geometric design parameters, operating conditions, and/or control inputs. Consequently, advances in our ability to interpret and predict the parameter-dependent dynamics of complex engineering systems will be a critical enabler for real-time control applications~\citep{brunton2015amr,herrmann2020se,herrmann2022prf}, many-query optimization tasks~\citep{gunzburger1999ijnmf,giles2000ftc,herrmann2016ate,herrmann2018hmt}, and the development of digital twins~\citep{grieves2017springer,hartmann2018springer,niederer2021natcs}.

Modal decomposition techniques are a cornerstone for the analysis of spatiotemporal phenomena, identifying dynamically relevant spatial patterns, or modes, that are associated with a simple time evolution. These analyses often enable building low-order representations of the underlying dynamics in terms of a combination of a few dominant modes, making it easier to understand and predict the future state of the system. In the context of structural dynamics, modal analysis refers to the solution of the eigenproblem arising from the linearized equations of motion to identify vibration modes and their associated natural frequencies~\citep{fu2001book}. Similarly, in thermal engineering, it is common practice to approximate the solution to transient heat conduction problems by considering only the slowest decaying temperature modes from the corresponding eigenproblem~\citep{bergman2011book}. Analogously, in the fluid dynamics community, the amplification, decay, and/or oscillation of perturbations to a given flow field is investigated by performing an eigendecomposition of the linearized dynamics operator using tools from stability theory~\citep{theofilis2011arfm}. Moreover, the unprecedented availability of high-fidelity data from numerical simulations, experimental measurements, and field recordings has led to the development and application of many data-driven modal decomposition methods to analyze spatiotemporal dynamical systems directly from data~\citep{taira2017aiaa,taira2019aiaa,herrmann2021jfm,schmid2022arfm,baddoo2023prsa}.

One prominent application of modal decomposition is in projection-based model order reduction~\citep{benner2015sr,rowley2017arfm}. In many engineering applications, computational models represent high-dimensional dynamical systems, thus involving a large number of degrees of freedom. The computational cost of these large models is often prohibitively expensive, making many-query tasks, such as design optimization or uncertainty quantification, where multiple instances of the same model need to be run, and real-time applications, such as feedback control and real-time digital twins, intractable. Model order reduction seeks to build reduced-order models (ROMs) that can be simulated at a fraction of the computational cost while approximately retaining the fidelity of the original full-order model~\citep{antoulasbook,bennerbook}. The eigendecomposition of a linear (or linearized) dynamical system offers one of the simplest approaches to model order reduction, known as modal truncation~\citep{antoulasbook}. This approach entails performing a Petrov-Galerkin projection of the governing equations onto a truncated set of eigenmodes, resulting in a diagonal operator for the reduced-order dynamics determined by the corresponding eigenvalues~\citep{antoulasbook,rowley2017arfm}.

In many scenarios, such as in design optimization or parameter estimation problems, it is necessary to simulate multiple realizations of a parameterized system. To approximate the system behavior across a range of parameters, parameterized reduced-order model have been introduced~\citep{benner2015sr}. One of the simplest approaches to create this is to compute a global basis of modes that encodes information of the entire parameter space~\citep{benner2015sr, legresley2000airfoil, weickum2009multi, amsallem2015design}. This has the drawback of requiring a large number of modes if the parameter space is large. Alternatively, local bases can be constructed for multiple sampled parameter values and then, either directly interpolated, or used to build multiple reduced-order operators that are themselves interpolated~\citep{benner2015sr}. Direct interpolation of the local bases involves treating the modes as continuous functions of the parameters. A more sophisticated approach, developed by~\cite{amsallem2008aiaa}, introduced mode interpolation on matrix manifolds, resulting in the preservation of mode properties, such as orthogonality. The latter method has been employed to interpolate reduced operators for the use of ROMs in parameterized systems, as well as to solve optimization problems~\citep{, degroote2010interpolation,amsallem2011online,amsallem2016real,choi2020gradient,hess2023data}. Interpolation on matrix manifolds interpolates each entry of the sampled matrix after a mapping to a matrix manifold. If the matrix to be interpolated is high-dimensional, the interpolation can quickly become computationally expensive. Additionally, it requires prior identification of a suitable matrix manifold onto which the matrix must be mapped and subsequently unmapped. For systems with controls, parameter-dependent dynamics are traditionally addressed using gain-scheduling, or with linear parameter varying (LPV) control to improve stability and robustness over the entire operating regime~\citep{toth2010modeling}. In a data-driven context, different approaches have been used, from learning low-order LPV models from a collection of linear time-invariant (LTI) systems~\citep{iannelli2021ast}, to building and simulating multiple ROMs and then interpolating their solutions in physical space~\cite{fonzi2020prsa}.

Variation of eigenmode shapes with parameters has been studied in fluid and structural mechanics~\citep{heinze2008robust, sorokin2015eigenfrequencies,xing2024dynamic}. In the context of structural mechanics, \emph{modal meta-modeling} has been developed to model variation of eigenvectors and eigenvalues with system parameters directly using regression analysis, and indirectly, solving a new eigenproblem instance~\citep{gallina2011enhanced,lu2019uncertainty,gibanica2021data,lee2024model}. In~\cite{gallina2011enhanced}, an indirect approach of meta-modeling is introduced, approximating each deformed eigenmode, avoiding the otherwise required mode pairing step. High dimensional linear and polynomial regression, as well as multi-output Gaussian processes have been used to model parameter-dependent modal data in structural systems in~\citep{gallina2011enhanced, lu2019uncertainty, gibanica2021data}. Dimensionality reduction using principal component analysis has been used in~\citep{gibanica2021data} to reduce the number of variables involved in the regression. In~\citep{lee2024model}, parameter-dependent eigenvectors are written using a low-rank approximation, enabling the solution of a reduced eigenproblem.

\begin{figure}[h!]
    \centering
    \includegraphics[width=17cm]{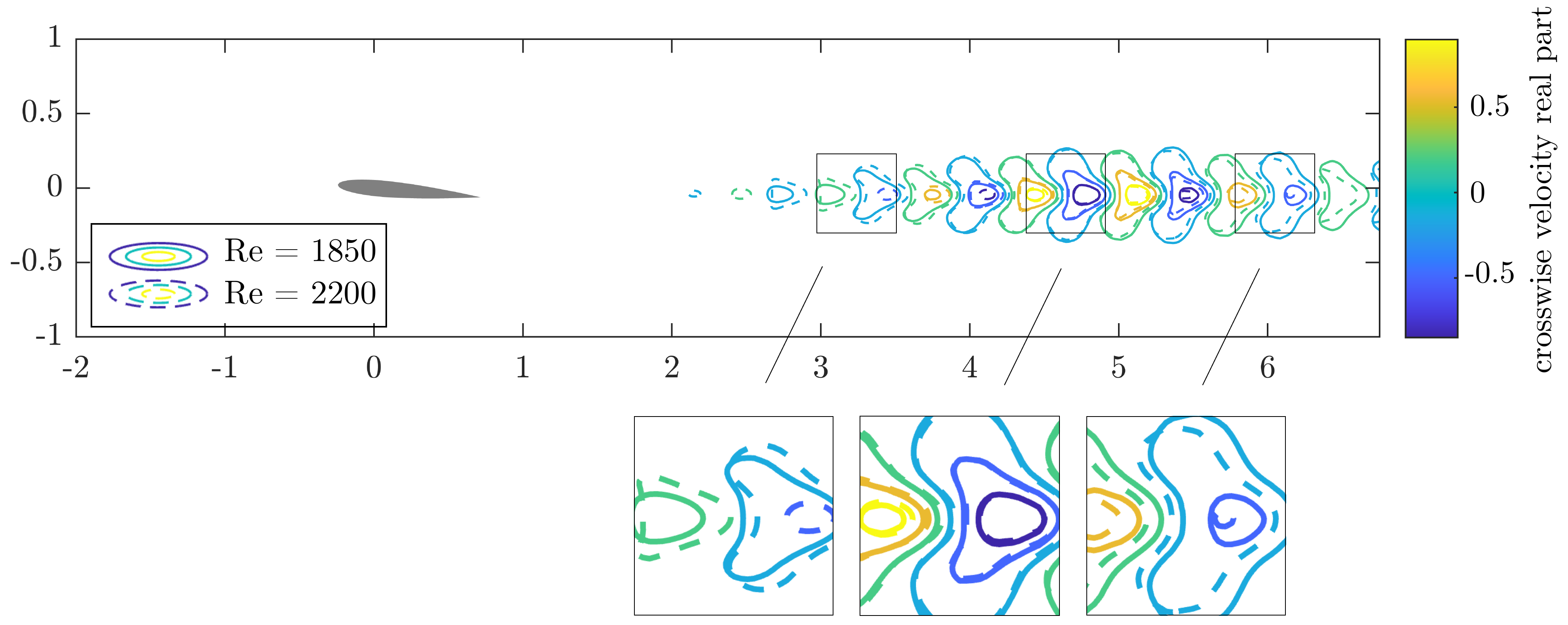}
    \caption{Crosswise velocity contours of the real part of the leading eigenmode for the incompressible flow over a NACA0012 airfoil. Contours of the eigenmodes at a Reynolds number $\mathrm{Re}=1850$ (solid lines) and $\mathrm{Re}=2200$ (dashed lines) are overlayed to show parameter-dependent mode deformation. Colors denote the same contour levels in both cases.}
    \label{ModalDeformation}
\end{figure}

Despite the numerous approaches for parametric model reduction, none of the aforementioned physics-based and data-driven modal analysis techniques incorporate parameter dependency in the extracted modes. For parameterized dynamical systems, the modes computed for one set of parameter values are not necessarily dynamically relevant for another set of parameters. Therefore, eigenmodes and eigenvalues extracted from modal analysis vary as the parameters change, which represents a challenge for the analysis of parameterized systems. This is exemplified for a fluid dynamics problem in Fig.~\ref{ModalDeformation}, showing the deformation of velocity eigenmodes. Furthermore, in the context of parameterized model reduction, capturing this mode deformation for a new parameter instance requires interpolating high-dimensional vectors, which quickly becomes computationally expensive if it has to be performed multiple times. Alternatively, if the reduced-order operators are interpolated instead, then the physical interpretability of the reduced state as a vector of modal amplitudes is lost.

\begin{figure}
    \centering
    \includegraphics[width=18cm]{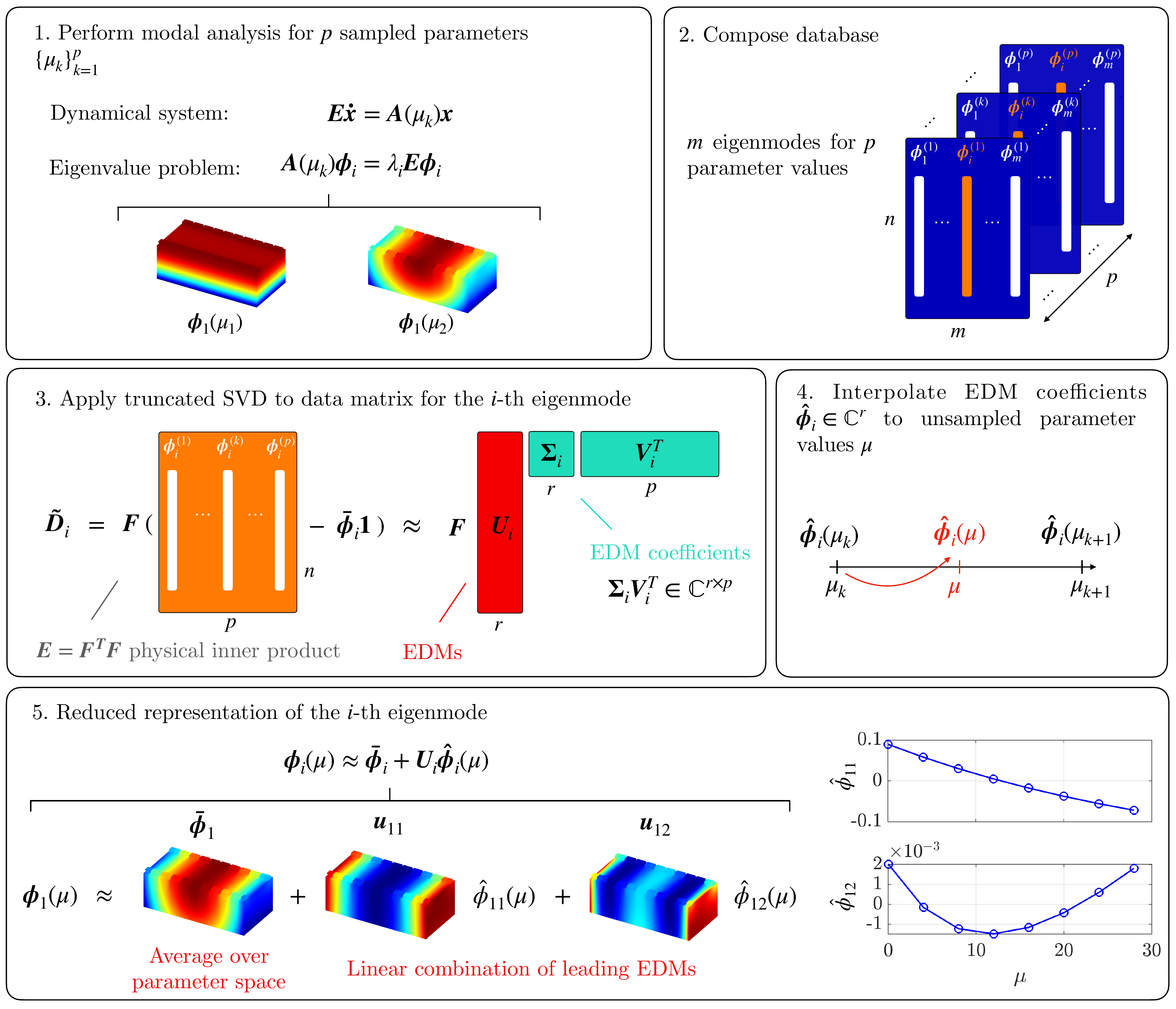}
    \caption{Schematic of the proposed method. An interpretable parameterized low-order representation of the eigenmodes of a system is built from a database of eigenmodes computed for a few sampled parameter values in the range of interest. As a result, we obtain eigen-deformation modes (EDMs) and EDM coefficients. In the example, the leading eigenmode of a heat conduction problem in a battery with parameterized boundary condition is expressed in a reduced form using two EDMs.}
    \label{SchematicMethod}
\end{figure}

In this work, we develop a method to produce an interpretable compressed representation of the deformation of eigenmodes with parameters in parameterized dynamical systems. To achieve this, the proposed method exploits correlations between eigenmodes computed for different parameter values to build a reduced-order representation of each eigenmode of interest in terms of a linear combination of eigen-deformation modes (EDMs). More specifically, EDMs are extracted from a singular value decomposition (SVD) to optimally (in a sense discussed in section~\S\ref{sec:methods}) capture the variation of eigenmodes across parameter space, as schematically shown in Fig.~\ref{SchematicMethod}. In addition, this enables building parameterized reduced-order models while avoiding expensive mode interpolation in physical space. 

The remainder of the paper is organized as follows. Theoretical background on modal analysis for parameterized dynamical systems and common approaches for parameterized projection-based model order reduction are covered in \S\ref{sec:background}. Our method to capture eigenmode deformation is formulated in~\S\ref{sec:methods}. The proposed method is demonstrated on three numerical examples in~\S\ref{sec:results} to show and discuss the resulting reduced-order eigenmode representations and their use in parameterized model reduction. Finally, we discuss our conclusions in~\S\ref{sec:conclusions}.

\section{\label{sec:background}Background}

\subsection{\label{sec:linearSystems}Modal analysis for parameterized systems}

Here we review the basics of modal analysis from the viewpoint of dynamical systems. Specifically, we consider parameterized linear dynamical systems of the form

\begin{equation}
    \b{E}\b{\dot{x}} = \b{A}(\mu) \b{x},
    \label{dynamicalSystem}
\end{equation}

where $\b{x}(t) \in \mathbb{R}^{n}$ is the state vector containing $n$ degrees of freedom, $t$ is time, the overdot denotes time differentiation, and $\b{E}, \ \b{A}\in \mathbb{R}^{n\times n}$ are the mass matrix and dynamics operator, respectively, with the latter incorporating the dependence on parameter $\mu\in\mathbb{R}$. Systems of ordinary differential equations like~\eqref{dynamicalSystem}, arise when modeling spatiotemporal phenomena governed by a partial differential equation (PDE), such as in heat transfer, structural mechanics, and fluid dynamics, and the spatial domain is discretized, for example using the finite element method. In those scenarios, the state vector $\b{x}$ is a discrete representation of the underlying continuous field. Importantly, the state here represents a perturbation from an equilibrium condition, therefore, as we show in appendix~\ref{app:therm}, Eq.~\eqref{dynamicalSystem} may also represent the dynamics of systems for which the underlying PDE has a source term, such as internal heat generation in thermal problems or a constant body force in fluids problems. Although for mechanical systems the equations of motion are of second order, they can also be rewritten in the form of~\eqref{dynamicalSystem} with a simple change of variables, as shown in appendix~\ref{app:mech}. We remark that the dynamics of the system will be linear either if the underlying PDE is linear, or if it has been linearized about a steady state equilibrium solution. Moreover, depending on the boundary conditions of the underlying PDE, the mass matrix $\b{E}$ may be singular (not invertible).




In eigenmode coordinates, the equations are decoupled, which allows building solutions using the principle of superposition. The time response of the system can be expressed as a linear combination of individual eigenmodes, that come from the eigenproblem

\begin{equation}
\label{eigenvalueProblem}
    \b{A}(\mu) \b{\phi}_i(\mu) = \lambda_i(\mu) \b{E}\b{\phi}_i(\mu),
\end{equation}

\noindent where $\b{\phi}_i(\mu) \in \mathbb{C}^{n}$, with $i=1,...,n$, is an eigenmode associated with the corresponding eigenvalue $\lambda_i(\mu) \in \mathbb{C}$. It is important to note that, in parameterized dynamical systems, eigenmodes and eigenvalues depend on the parameters due to the parameter dependence of the system matrices, in this case only considered in $ \b{A}(\mu)$. The time evolution of the system, in physical coordinates, is given by

\begin{equation}\label{spectralSolution}
    \b{x}(t, \mu) = \b{\Phi}(\mu) \exp{(\b{\Lambda }(\mu)t})\b{\Phi}(\mu)^{-1} \b{x_0},
\end{equation}

\noindent where $\b{\Phi}(\mu) \in \mathbb{C}^{n \times n}$ is a matrix containing the eigenmodes $\b{\phi}_i(\mu)$ in its columns, $\b{\Lambda}(\mu) \in \mathbb{C}^{n \times n}$ is a diagonal matrix containing the eigenvalues corresponding to each column of $\b{\Phi}(\mu)$ in its diagonal. Equation (\ref{spectralSolution}) represents the exact solution of the system with parameter $\mu$ from an initial condition $\b{x_0} \in \mathbb{R}^{n}$ as a combination of the eigenmode exponential dynamics.

\subsection{\label{sec:eigenmodeTruncation}Eigenmodes in model order reduction}

An approximation of system~\eqref{dynamicalSystem} can be formulated if equations are projected onto a lower-dimensional subspace by employing the so-called reduced-order bases, constructing a projection-based ROM~\citep{antoulasbook, rbbook, loewnerbook, bennerbook, benner2015sr, rowley2017arfm}. In projection-based reduced-order modeling, a common choice of basis is a set of $m$ direct eigenvectors and adjoint eigenvectors arising from the associated eigenproblem of the system, which we denote as $\b{\Phi_m}(\mu)$ and $\b{\Psi_m}(\mu) \in \mathbb{C}^{n\times m}$, respectively. Notice that we added the parameter dependence in the reduced-order bases, due to the dependence of the eigevectors on the parameter. By the above election, the solution is being approximated by considering the first $m$ dominant eigenmodes. Approximating the full-dimensional state as an evolution on the $m$ dimensional subspace spanned by the retained eigenvectors, we write  $\b{x} = \b{\Phi_m}(\mu)\b{\hat x}$, where $\b{\hat x} \in \mathbb{C}^{m}$ is the reduced state. In this case, reduced bases are bi-orthogonal, and if properly normalized, they satisfy $\b{\Psi_m}^{^{\mathrm{H}}}\b{E\Phi_m} = \b{I}$, where $\b{I}\in\mathbb{R}^{m\times m}$ is an identity matrix and $\mathrm{H}$ denotes the Hermitian transpose. Consequently, the reduced system takes the simple form

\begin{equation}\label{DiagonalizedSystem}
     \dot{\b{\hat{x}}} = \b{{\Lambda}_m}(\mu) \b{\hat x},
\end{equation}

\noindent where $\b{{\Lambda}_m}(\mu) \in \mathbb{C}^{m\times m}$ is the diagonal matrix containing the $m$ eigenvalues corresponding to the columns of $\b{\Phi_m}$ in its diagonal. In this case, the approximation of the high-dimensional system based on this $m$-dimensional subspace projection is

\begin{equation}\label{ROMApproximation}
    \b{x}(t, \mu) \approx \b{\Phi_m}(\mu) \exp{(\b{{\Lambda}_m}(\mu)t}) \b{\hat x_0},
\end{equation}

\noindent where $\b{\hat x_0}$ is the initial condition in reduced coordinates, given by $\b{\hat x_0} = \b{\Psi}^{\mathrm{H}}(\mu) \b{E x_0} $. The ROM under these assumptions can be simulated from a given initial condition if the bases of eigenmodes, $\b{\Phi_m}(\mu)$ and $\b{\Psi_m}(\mu)$, and eigenvalues in $\b{\Lambda_m}(\mu)$ are known at the parameter $\mu$. The above represents a specific application where parameterized eigenmodes would be helpful. Furthermore, a reduced representation of eigenmodes would be especially valuable in cases where the ROM must be computed multiple times for different parameter values, such as in design optimization problems~\citep{benner2015sr}.

\section{\label{sec:methods}Proposed method}

In this section we develop our method to build a parameterized reduced representation of the eigenmodes of a system using modal analysis data. 

\subsection{\label{subs:MethodDetails}Optimal basis for eigenmode deformation}

Consider a linear dynamical system of the form of~\eqref{dynamicalSystem}, that can arise from the linearization of a non-linear system around an equilibrium point, parameterized by parameter $\mu \in \Omega_\mu$ and $n$ state variables. The proposed method takes the case of a single parameter, so $\Omega_\mu \subset \mathbb{R}$ is an interval. Let $\{ \mu_k \}_{k=1}^p$ denote a set of $p$ sampled parameters in $\Omega_\mu$. Here we assume equi-spaced in $\mu$. Adapting the proposed method to non-uniformly sampled parameters is straightforward but is not addressed in this work to maintain simplicity. We use $\b{\phi}_i^{(k)}$ to denote the $i$-th eigenmode evaluated at parameter $\mu_k$. This process requires grouping eigenmodes into sets according to a physically consistent order for different parameters. The aforementioned \emph{mode pairing} usually can be done based on eigenvalue ordering. However, special care must be taken in cases of mode degeneration \cite{gallina2011enhanced}. The above is discussed in more detail in section~\S\ref{sec:Interpolation}.

The proposed method utilizes a database comprised of the first
$m$ eigenmodes computed for each sampled parameter $\mu_k$ (correctly paired), and the non-parameter dependent mass matrix $\b{E} \in \mathbb{R}^{n \times n}$ (if it appears in the analyzed dynamical system). 
%
%
%
%
For each $i$-th eigenmode, with $i=1,...,m$, the data matrix $\b{D}_i \in \mathbb{C}^{n\times p}$ is constructed as

\begin{equation}\label{VjdefMatrix}
\b{D}_i = 
\begin{bmatrix}
\vert & \vert & & \vert \\
\b{\phi}_i^{(1)}   & \b{\phi}_i^{(2)} & \hdots & \b{\phi}_i^{(p)} \\
\vert & \vert & & \vert
\end{bmatrix} - \b{\bar \phi}_i \b{1}^{T},
\end{equation}

\noindent where $\b{\bar \phi}_i$ denotes the mean $i$-th eigenmode (across parameter samples) $\b{\bar \phi}_i = \frac{1}{p} \sum_{k=1}^p \b{\phi}_i^{(k)}$, and $\b{1} \in \mathbb{R}^n$ is a constant vector of ones. Each column of the matrix $\b{D}_i\in \mathbb{C}^{n \times p}$ is the difference between the $i$-th eigenmode for a sampled parameter $\mu_{\mathrm{k}}$ and the mean eigenmode $\b{\bar \phi}_i$, so $\b{D}_i$ contains the variation of the $i$-th eigenmode fluctuation from the mean over the parameter interval $\Omega_\mu$.

To ensure the use of a spatial inner product that is physically meaningful during the underlying optimization problem—accounting for uneven spatial discretizations and potentially different quantities in the state vector—the modified data matrix $\b{\tilde{D}}_i$ is computed as follows:

\begin{equation}
    \b{\tilde{D}}_i =  \b{FD}_i
\end{equation}

\noindent where $\b{F} \in \mathbb{R}^{n \times n}$ is the Cholesky factor of the mass matrix $\b{E}$, so it satisfies $\b{F}^T\b{F}=\b{E}$.

Singular value decomposition (SVD) is then applied to $\b{\tilde{D}}_i$, retaining $r$ singular vectors and singular values, resulting in a rank-$r$ optimal approximation of $\b{\tilde{ D}}_i$, in the sense of the Frobenius norm \cite{databook},

\begin{equation}
    \b{\tilde{D}}_i \approx \b{\tilde U}_i \b{\Sigma}_i \b{V}_i^{T}.
\end{equation}

\noindent Various criteria exist to select $r$ in SVD-based dimensionality reduction, such as capturing a specific percentage of variance or applying the optimal hard thresholding criterion from~\citep{gavish2014optimal}.

The $i$-th set of \emph{eigen-deformation modes} (EDMs) are obtained as

\begin{equation}
    \b{U}_i = \b{F}^{-1} \b{\tilde U}_i,
\end{equation}

\noindent and the matrix $\b{\Sigma}_i \b{V}_i^{T} \in \mathbb{C}^{r \times p}$ contains the evolution of EDM coefficients $\b{\hat \phi }_i^{(k)} \in \mathbb{C}^{r}$ over the sampled parameter interval $\Omega_\mu$, so the original $\b{D}_i$ data matrix can be approximated as

\begin{equation}\label{LowRankAprox}
    \b{D}_i \approx \b{U}_i \b{\Sigma}_i \b{V}_i^{T}.
\end{equation}

The matrices $\b{U}_i \in \mathbb{C}^{n \times r}$, for $i=1,...,m$, contain $r$ EDMs $\b{u}_{ij}$ for $j=1,...,r$, in their columns that form an orthogonal base that optimally captures the deformation of the $i$-th eigenmode, in the sense of a physically meaningful spatial inner product averaged across the parameter interval. Similarly to eigenmodes, EDMs are spatial patterns in state space, so their shape can be visualized and analyzed. This formulation allows approximating each eigenmode using $r$ EDMs along with their corresponding EDM coefficients

\begin{equation}
    \b{\phi}_i(\mu_k) \approx \b{\bar \phi}_i + \b{U}_i \b{\hat{\phi}}_i^{(k)}.
\end{equation}

The mode shapes for an unsampled parameter $\mu$ can be obtained by interpolating the coefficients vector $\b{\hat{\phi}}_i^{(k)}$, leading to a parameterized low-order representation of each eigenmode:

\begin{equation}\label{ReducedRepresentation}
    \b{\phi}_i(\mu) \approx \b{\bar \phi}_i + \b{U}_i \b{\hat \phi}_i(\mu).
\end{equation}

\subsection{\label{RelationPOD}Relation to proper orthogonal decomposition}

When formulated in terms of continuous functions, the proposed method provides a basis for the modal decomposition of an ensemble of functions, and so it shares a tight connection with the proper orthogonal decomposition (POD)~\citep{berkooz1993arfm}. Let $\phi_i(x,\mu)$ be the $i$-th eigenfunction of the linear (or linearized) underlying PDE that upon discretization leads to Eq.~\eqref{dynamicalSystem}, which depends on a spatial variable $x$ and the parameter $\mu$. Given an ensemble of instances of $\phi_i$ obtained for various parameter values, we ask what is the function $u_i(x)$ that is most similar to the members of our ensemble on average. To provide a satisfactory answer we need to mathematically define a notion of closeness between functions and an averaging operation. The former is usually assessed in terms of alignment based on the spatial inner product
\begin{equation}
    \left(f(x),g(x) \right) = \int_{\Omega_x} \! f(x)g^*(x)\mathrm{d}x,
\end{equation}
where $f$ and $g$ are square integrable functions over the spatial domain $\Omega_x$ and the asterisk denotes complex conjugation. Moreover, in the context of this work, we use the following averaging operation over the parameter dependence
\begin{equation}
    \langle h(\mu)\rangle = \frac{1}{|\Omega_{\mu}|}\int_{\Omega_{\mu}} \! h(\mu)\mathrm{d}\mu,
\end{equation}
where $h$ is a square integrable function over the parameter domain of interest $\Omega_{\mu}$, and $|\Omega_{\mu}|=\int_{\Omega_{\mu}} \! \mathrm{d}\mu$ is the length of the interval.

With these inner product and averaging definitions, the search for $u_i(x)$ can be formulated as the following optimization problem
\begin{equation}
    u_i = \underset{u_i'}{\mathrm{argmax}}\frac{\left\langle |\left(\phi_i,u_i'  \right)|^2 \right\rangle}{\left( u_i', u_i' \right)},\label{eq:pod_opt}
\end{equation}
which is equivalent to the POD problem~\citep{berkooz1993arfm}, except our ensemble is formed by observations of the eigenfunctions $\phi_i$ over different parameter instances. As for POD~\citep{berkooz1993arfm}, the optimal and suboptimal solutions to Eq.~\eqref{eq:pod_opt} are given by the eigenfunctions of the following Fredholm integral eigenvalue problem
\begin{equation}
    \int_{\Omega_x} \! \left\langle \phi_i(x) \phi_i^*(x')\right\rangle u_{ij}(x') \mathrm{d}x = \sigma_{ij}^2 u_{ij}(x),\label{eq:pod_eig}
\end{equation}
where $\sigma_{ij}^2$ and $u_{ij}$ correspond to the $j$-th largest eigenvalue and the corresponding eigenfunction of the two-point spatial correlation tensor $R_i(x,x')=\left\langle \phi_i(x) \phi_i^*(x')\right\rangle$ for our ensemble of the $i$-th eigenfunction of the system. Consequently, the $u_{ij}$ functions — that are the continuous formulation of EDMs — form an optimal basis to represent $\phi_i$ by exploiting spatial correlations across the parameter interval of interest.

\section{\label{sec:examples}Numerical examples and dataset}

\begin{figure}
    \centering
    \includegraphics[width=17cm]{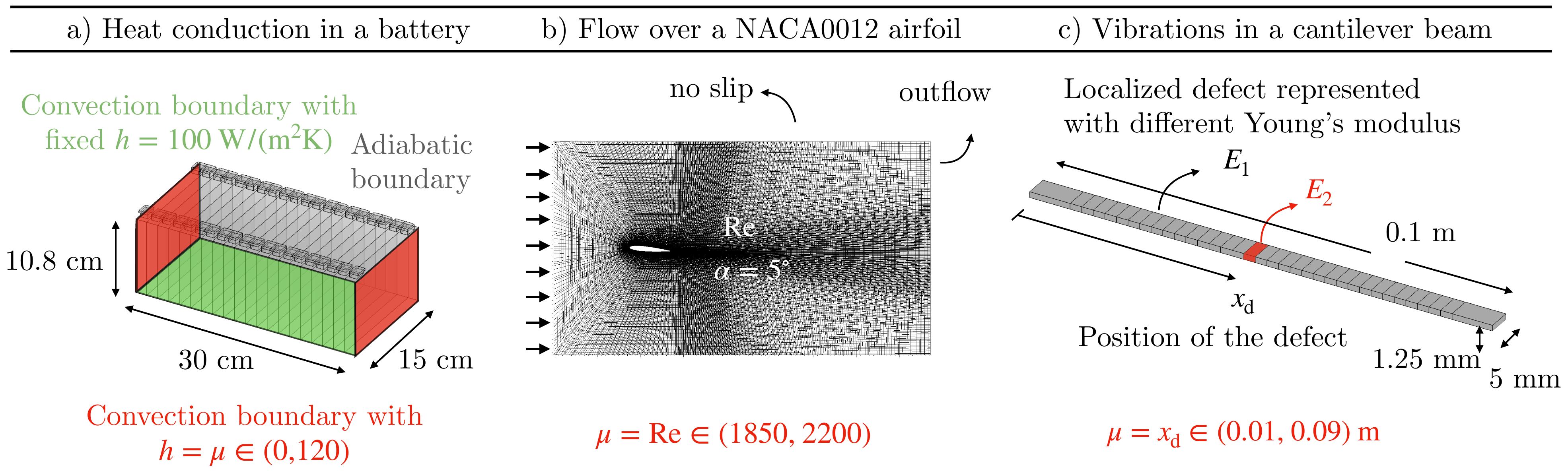}
    \caption{Schematic the numerical examples used to test the proposed method. (\textit{a}) Heat conduction in a battery parameterized by the convective heat transfer coefficient acting on the back and front faces shown in red. (\textit{b}) Flow over a NACA0012 airfoil parameterized by the Reynolds number.  (\textit{c}) Vibrations in a cantilever beam parameterized by the longitudinal position of a localized defect.}
    \label{ExamplesSchematics}
\end{figure}

To demonstrate the application of the proposed method, we generate datasets comprised of eigenmodes computed for different parameters of three systems governed by parameterized PDEs commonly found in engineering: the heat conduction in a modular battery, the flow over a NACA0012 airfoil, and the mechanical vibrations in a cantilever beam. Schematics of the three examples are shown in Fig.~\ref{ExamplesSchematics}, and all the data is available on github.com/ben-herrmann. The modes for different parameters are paired based on the (descending) order of the real part of the eigenvalues. In the fluid mechanics example, however, special measures are required to address the observed eigenvalue crossing~\citep{gallina2011enhanced}, further elaborated in section~\S\ref{sec:Interpolation}. All the eigenmodes extracted for these systems are pre-processed, as shown and discussed in appendix~\ref{app:preProcessing}, in order to set a consistent normalization over the sampled parameter interval.

\subsection{\label{heatTransferExample}Heat conduction in a battery}
For our first example, we consider the heat conduction in a battery with internal heat generation. Two types of boundary conditions are considered on the battery. Convection on the bottom, back, and front faces, and insulation (adiabatic boundary) on the rest of the domain. The battery is cooled from the bottom, with a fixed convective coefficient $h =\: \si{100\: \si{\W \per (m ^{2}K})}$. The cooling of the back and front faces occurs with a different convective coefficient $\mu$, that is the parameter considered for this example. The schematics of the boundary conditions is shown in Fig.~\ref{ExamplesSchematics} (\textit{a}). The ambient temperature is constant and equal to $T_{\infty}=293\si{\K}$.

The three-dimensional heat equation along with the respective boundary conditions is discretized via finite elements to obtain a system of differential equations of the form of~\eqref{dynamicalSystem}. Modeling, discretization, and modal decomposition is performed using MATLAB's PDE toolbox. The geometry of the battery is created using a function to model a modular prismatic battery, included in the toolbox. The details on the fixed geometric parameters and material properties used are shown in appendix~\ref{app:BatteryParameters}. The spatial domain is discretized using the default settings of the MATLAB mesh generating function, resulting in a mesh of $n = 28673$ tetrahedral elements. The dataset for this example is comprised of the first $10$ eigenmodes computed for parameters ranging from $\mu=0$ to $\mu=28$ in 100 equal increments, adding the values for $\mu=20,40,80,120$. 

\subsection{\label{fluidMechanicsExample}Flow over a NACA0012 airfoil}
In our second example, we consider the two-dimensional flow over a NACA0012 airfoil governed by the incompressible Navier-Stokes equations. We choose a fixed angle of attack of $\alpha=5^{\circ}$ and use the Reynolds number as our parameter. The Reynolds number is a dimensionless quantity used in fluid mechanics that compares the relative importance of inertial and viscous effects and, for the flow over an airfoil, is defined as $\mathrm{Re}=U c/\nu$, with $U$, $c$ and $\nu$ being the free-stream velocity, the airfoil chord and the kinematic viscosity, respectively. We restrict our attention to the range $\mathrm{Re}\in(1850,2200)$, where the flow is laminar and at $\mathrm{Re}=2130$ undergoes a bifurcation from a stable equilibrium (steady flow) to an unstable equilibrium that gives rise to an unsteady flow characterized by self-sustained periodic vortex shedding~\cite{gupta2023jfm}. The eigenmodes whose parameter dependence we seek to compress arise from the linearization of the incompressible, two-dimensional Navier–Stokes equations about the (stable and then unstable) equilibrium flow. These eigenmodes represent two-dimensional perturbations (constant in the spanwise direction) to that same equilibrium velocity field. 

To build a dataset of eigenmodes for various parameter values, we need to compute the equilibrium flow and subsequently perform a stability analysis for every parameter value. For this purpose, we use the spectral element code Nek5000~\cite{fischer2008nek5000} in conjunction to the recently developed nekStab toolbox~\cite{frantz2023amr}. We consider a rectangular computational domain extending from $-2c$ to $6c$ in the streamwise direction and $-2.5c$ to $2.5c$ in the crosswise direction, with the leading edge of the airfoil located at the origin. Spatial discretization relies on a C-grid mesh embedded inside the rectangular domain, as shown in Fig.~\ref{ExamplesSchematics}(b), using 4368 spectral elements with a polynomial order of $N = 5$. The stable equilibria for values of $\mu=\mathrm{Re}<2130$ are computed by running direct numerical simulations (DNS) starting from rest over a time horizon of $500$ time units. The unstable equilibria, for $\mu=\mathrm{Re}\ge2130$, are computed using the selective frequency damping (SFD) method~\cite{aakervik2006pof} implemented in nekStab using a filter gain of $\xi=0.05$ and a cutoff frequency of $\omega_c=0.12$. Specifically, a DNS is started from rest and evolved over $200$ time units and, from that state, SFD is used to simulate an additional $200$ time units. The leading $m=40$ eigenmodes and eigenvalues are computed using an Arnoldi iteration, also implemented in nekStab~\cite{frantz2023amr}, with a Krylov basis dimension $m=350$ and a time step $\tau=0.1$. The dataset generated for this example includes the equilibria, and $40$ eigenmodes and eigenvalues for each of eight parameter values that go from $\mu=1850$ to $2200$ in increments of $50$. Vectors containing the grid coordinates and the cell volumes are also included in the dataset to be used for plotting and as integration quadrature weights, respectively.

\subsection{\label{structuralExample}Vibrations in a cantilever beam}
For the last example, we consider the vibrations in a cantilever beam that has a localized defect, with the position of the defect as a parameter. This system is governed by the linear elasticity equations, and the localized defect is represented as a region with a different Young's modulus than the rest of the beam. Discretization is applied via finite elements to obtain a system of the form of~\eqref{dynamicalSystem}.

The beam has dimensions of $10 \si{\cm} \times 5 \si{\mm} \times 1.25 \si{\mm}$ and is made of steel with Young's modulus $E_1 = 210\:\si{\giga \pascal}$, mass density $\rho =\: \si{7800\: \si{\kilogram \per m ^{3}}}$ and Poisson ratio $\nu = 0.3$. The parameter of this example is the longitudinal position of the defect on the beam, in meter units. The defect is modeled as a cross-sectional region, $1.25 \si{\mm}$ in length, with a different Young's modulus, $E_2 = 10\:\si{\giga \pascal}$, located at a specific position along the $x$-axis. The schematics of the beam is shown in Fig.~\ref{ExamplesSchematics} (\textit{c}). Modeling, discretization, and modal analysis are performed using MATLAB's PDE toolbox. The geometry is created with 65 cells, each $1.25 \si{\mm}$ long, located at positions ranging from $0.01\: \si{\m}$ to  $0.09\: \si{\m}$ along the $x$-axis, denoted as $x_\mathrm{d}$. One of the aforementioned cells is assigned to have a different Young's modulus than the rest of the beam before solving a modal analysis, allowing for data collection. The spatial domain is discretized using $10375$ tetrahedral elements using MATLAB's mesh generating function, setting a maximum edge length value of 0.002. The data set is comprised of the first $5$ eigenmodes for $65$ parameter values ranging from $\mu=0.01$ to $\mu=0.09$ in increments of $1.25 \cdot 10^{-3}$.

\section{\label{sec:results}Results and discussion}

In this section, we demonstrate the application of the proposed method using the datasets described in the previous section. First, we apply the method to show the resulting EDMs and EDM coefficients for the first eigenmode of each system. We then investigate the effectiveness of the EDMs in compressing parameter-dependent eigenmode deformation, comparing our reduced representation to the interpolation of eigenmodes in physical space. Finally, we demonstrate the use of EDMs for parameterized model reduction, using the heat transfer example in a limited data scenario, shown in Fig. \ref{ROMResultsFigure}.

\subsection{\label{sec:ModesResults}EDMs and EDM coefficients}

\begin{figure}
    \centering
    \includegraphics[width=17cm]{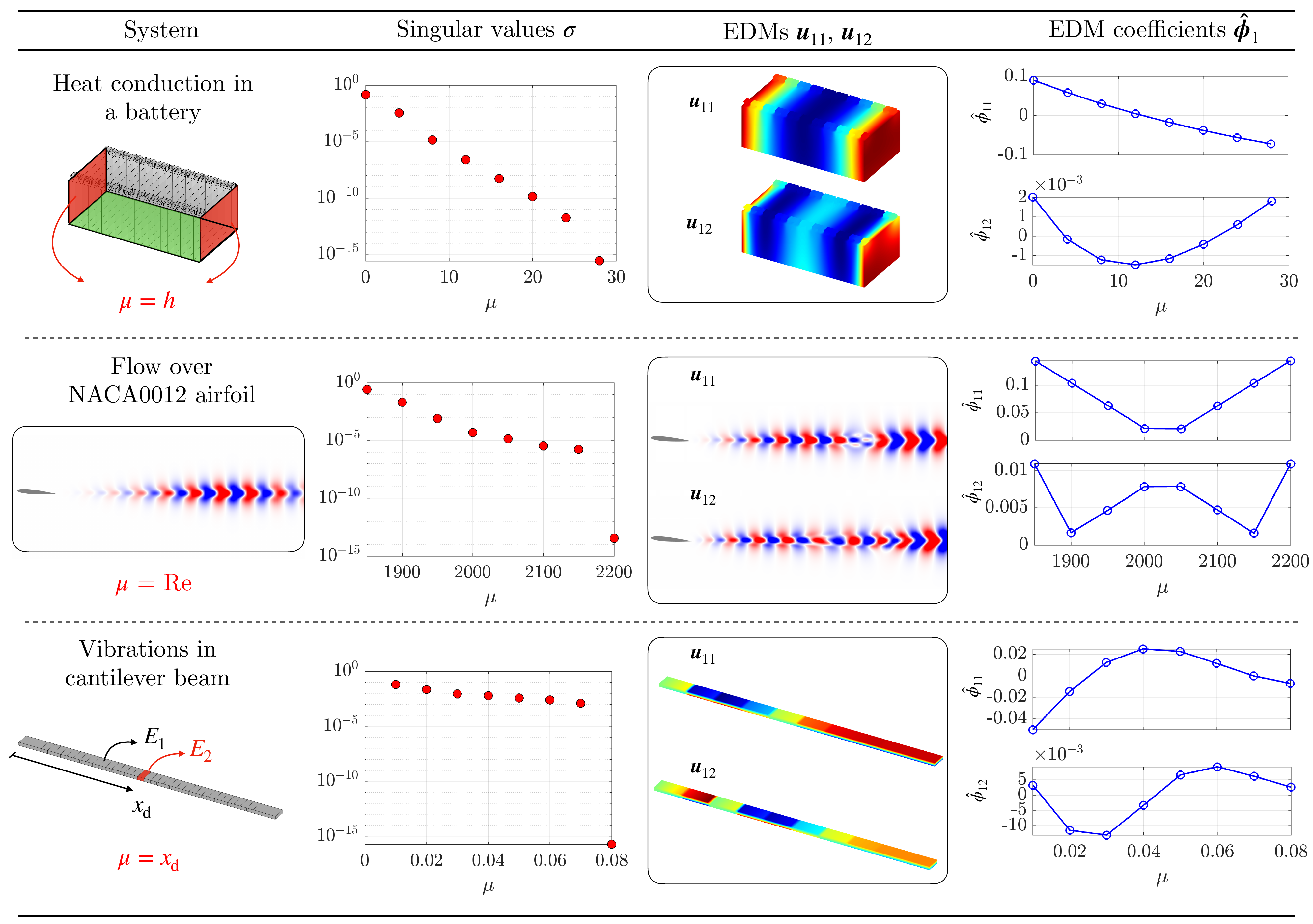}
    \caption{Application of the proposed framework to datasets for the leading eigenmode in parameterized heat transfer, fluid mechanics, and structural dynamics problems. The decay of the singular values of the eigenmode data matrix, the leading two EDMs, and their corresponding EDM coefficients are shown in the second, third, and fourth columns, respectively. For the fluids example, the magnitude of the complex-valued modes and coefficents are displayed.}
    \label{ModesResultsFigure}
\end{figure}

The results presented in this section are computed using eight parameter values for each example. For the battery case, the parameters range from $\mu = 0$ to $\mu = 28$ in increments of 4. For the airfoil, the parameters range from $\mu = 1850$ to $\mu = 2200$ in increments of $50$. Lastly, for the cantilever beam, the parameters range from $\mu = 0.01$ to $\mu = 0.08$ in increments of $0.01$. The results of applying the proposed method to each data set are shown in Fig.~\ref{ModesResultsFigure}. The singular values $\sigma$, resultant EDMs, and their corresponding coefficient evolution over the sampled parameter interval for the first eigenmode are shown. EDMs are visualized as temperature and $x$-displacement fields for the heat transfer and structural mechanics problems, respectively. For the fluid mechanics problem, the cross-wise velocity component is plotted, showing its real part. The fraction of energy captured from the data by retaining $r$ EDMs is defined by

\begin{equation}
    \% \mathrm{Energy} = \frac{\sum\limits_{k=1}^r \sigma_k}{\sum\limits_{k=1}^p \sigma_k}.
\end{equation}

In the case of the heat transfer example, $97.6\%$ is captured by retaining just one EDM. Retaining the first two EDMs, $99.9\%$ is captured, a considerable amount of energy to provide an accurate low-rank parameterized approximation of the first eigenmode.  In this example, when $\mu=0$, the first eigenmode is a vertically stratified temperature pattern because the back and front faces of the battery are insulated, so the battery is just being cooled from the bottom. As $\mu$ increases, convection starts to play a role in the direction perpendicular to the parameter faces, so the first eigenmode deforms, incorporating a gradient in that direction. The first and second EDMs effectively capture the deformation of this eigenmode from the average shape to another parameter within the sampled interval. Looking at the first EDM coefficient, we can see that it is approximately zero in the middle of the sampled interval ($\mu=15$), which indicates that the eigenmode at $\mu=15$ is roughly equal to the average sampled eigenmode. It is important to note that in this example, the parameter is a boundary condition, so the EDMs are essentially shifting the problem from one boundary condition to another. The first EDM represents the dominant effect of the boundary condition in the system, which is a temperature pattern that varies in the direction perpendicular to the convective boundaries involved in the parameter change. The EDMs effectively capture the deformation of the dominant dynamic pattern, allowing the transition of this pattern from one boundary condition to another through the addition of a linear combination of them.

In the case of the fluid mechanics example, the variance captured for the first eigenmode is $92.7\%$, $99.7\%$ and $99.9\%$ for values of $r$ of $1,2$ and $3$, respectively. The first mode exhibits structures that grow or shrink depending on the spatial region as the Reynolds number varies. Specifically, as shown in Fig.~\ref{ModalDeformation}, the structures increase in size and become concentrated near the airfoil, while they appear to decrease in size and become more diffuse farther from the airfoil. Looking at the obtained EDMs, we can observe similar structures to the eigenmodes both in the region near the airfoil and in the region far from the airfoil. Additionally, a more complex pattern emerges near the central area along the $x$-axis of the discretized domain. This indicates that, in this central region, the first eigenmode exhibits a more intricate deformation than the extreme regions mentioned earlier. The EDMs successfully capture a deformation that is not easily noticeable at first glance in this area, providing a more accurate parametric representation than simply using the fixed average eigenmode. By analyzing the deformation mode coefficients, we can see that the first coefficient is non-zero across the sampled interval, indicating that it is always necessary to add some deformation to the average eigenmode to properly represent the first eigenmode parametric variation.

Finally, in the case of the structural mechanics example, the first eigenmode deforms along the main axis of the beam as the localized defect moves along it. The decay of the singular values is significantly slower than that of their counterparts. The energy captured is $59.3\%$, $79.8\%$, $88.0\%$, $93.3\%$, $96.6\%$  and $98.8\%$ for values of $r$ of $1,2,3,4,5$ and $6$, respectively. Compared to the previous cases, in this type of example, a larger number of EDMs is required to provide an accurate reduced representation of the leading eigenmode. This occurs because the eigenmode in question exhibits a spatial structure that shifts as the parameter changes, behaving like a traveling wave (through parameter space), and our method relies on linear dimensionality reduction through an SVD, which is known to struggle with traveling waves~\cite{databook}. More generally, the compression capabilities of the proposed approach are limited by the rank of the underlying two-point spatial correlation tensor defined with the averaging operation across the parameter interval. Future research could investigate the use of nonlinear dimensionality reduction to improve the compression in these scenarios, albeit at the cost of loosing the interpretability of EDMs.

\subsection{\label{sec:Interpolation}EDM-based interpolation}

Here, we provide a more quantitative assessment of the accuracy of the reduced representation of the eigenmodes provided by our method. Specifically, eigenmode interpolation based on EDMs, given by Eq.~\eqref{ReducedRepresentation}, is compared to the direct interpolation in physical space of the eigenmodes in the three examples for different values of the number of retained EDMs $r$. The same eight sampled parameters as those used in the previous subsection are used for the battery example. For the cantilever beam, nine sampled parameters are used from $\mu=0.01$ to $\mu=0.09$ in increments of $0.01$. Finally, for the airfoil example, four sampled parameters are used, from $\mu=1850$ to $\mu=2150$ in increments of $100$. Performance is quantified in terms of the interpolation error defined as

\begin{equation}
    \mathrm{error} = \frac{\| \b{F}(\b{\phi}_i(\mu) - \b{\tilde \phi}_i(\mu))\|_2}{\| \b{F}\b{\phi}_i(\mu)\|_2},
\end{equation}

\noindent where $\b{\phi}_i(\mu)$ and $\b{\tilde \phi}_i(\mu)$ are the $i$-th ground truth eigenmode and its interpolated prediction computed at the unsampled parameter $\mu$, respectively. As earlier, $\b{F}$ is the Cholesky factor of the mass matrix used to maintain a physically meaningful inner product.

\begin{figure}
    \centering
    \includegraphics[width=17cm]{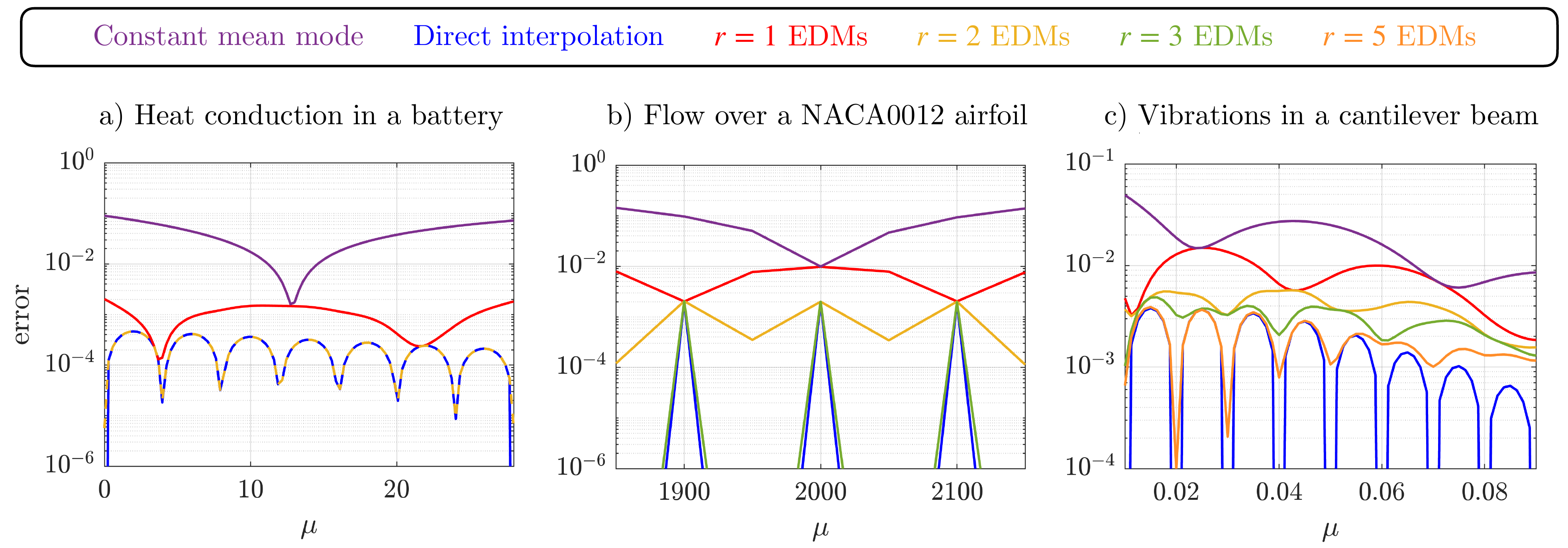}
    \caption{Spatially integrated error between the leading eigenmode and its interpolation over parameter space from a few sampled parameter values for the (\textit{a}) heat transfer, (\textit{b}) fluid mechanics, and (\textit{c}) structural dynamics examples. Direct interpolation of the eigenmodes in physical space is compared to EDM-based interpolation using a different number of EDMs.}
    \label{InterpolationErrorFigure}
\end{figure}

For the heat transfer example, the validation ground truth eigenmodes are computed in the dataset for $100$ parameter values covering the range of the training data, that is, from $\mu = 0$ to $\mu = 28$. For the structural mechanics example, $65$ positions of the defect from $\mu = 0.01$ to $\mu = 0.09$ detailed in section~\ref{sec:examples} are used as ground truth parameters to calculate the error. Finally, for the fluid mechanics example, seven parameters were available in the range of the training data, from $\mu=1850$ to $\mu=2150$ in increments of 50. Direct and EDM-based interpolations were performed using a linear interpolation scheme. For the heat transfer example, our reduced representation can match the error of the direct interpolation by retaining $r=2$ EDMs for the leading eigenmode, as shown in Fig.~\ref{InterpolationErrorFigure} (\textit{a}). Furthermore, this representation also outperforms the model using $r=1$ EDMs and the constant mean mode over the parameter interval ($r=0$ EDMs), as expected. For the airfoil example, a value of $r=3$ is sufficient to match the direct interpolation, as shown in Fig.~\ref{InterpolationErrorFigure} (\textit{b}). The error curve is similar to that of the battery example, but less smooth due to the limited number of validation parameters available. For the structural dynamics example, a bigger value of $r$ is needed to match the directly interpolated eigenmodes, as shown in Fig.~\ref{InterpolationErrorFigure} (\textit{c}). This is expected from the slow decay of the singular values, presented in the previous subsection, that arises due to the type of mode deformation involved in this problem. The same trends are observed for the third eigenmode of each system, and those results are included in appendix~\ref{app:ThirdEigenmode}. Our reduced representation converges to the direct interpolation for a sufficiently large value of $r$. EDMs effectively compress parametric mode deformation in all three engineering examples, enabling a reduced representation of the eigenmodes that is as accurate as the full-dimensional interpolation in physical space.

It is important to recall that a robust mode pairing strategy is critical to adequately build the eigenmode data matrices for the application of the proposed method. In many scenarios, mode pairing based on eigenvalue ordering, for example descending according to their real part, is enough to correctly follow specific eigenmodes through parameter space. However, in certain cases where mode degeneration~\cite{gallina2011enhanced} occurs, relying solely on eigenvalue ordering will lead to an incorrect mode pairing. We encounter this behavior for our fluid mechanics example, where eigenvalues cross their paths in the complex plane as the Reynolds number changes, resulting in the associated eigenmodes switching their order. This is easily addressed by pairing eigenmodes with similar associated eigenvalues and manually verifying mode shape similarity. No mode crossing is observed in the other two examples and mode pairing based on eigenvalue ordering is used. Future research will investigate the use of more sophisticated mode pairing strategies, such as mode shape similarity or mode tracking~\citep{gallina2011enhanced, lu2020mode} to deal with potential edge cases. 

\subsection{\label{sec:ROMResults}Parameterized reduced-order modeling using EDMs}
We use the reduced representation of the eigenmodes based on EDMs to build parameterized ROMs and compare their performance with other simple alternatives commonly used. The heat transfer problem is used as a test bed for this purpose, simulating the temperature evolution inside the battery. Four parameters were selected for training, namely, $\mu_k = 0, 40, 80, 120$. 
We consider a database for ROMs simulation that includes the sampled parameters $\mu_k$, the corresponding equilibria $\b{\bar x} (\mu_k)$ used as linearization points, and $m=6$ eigenmodes and eigenvalues for each sampled parameter to apply the model reduction.


ROMs are built using Eq.~\eqref{ROMApproximation} using interpolated eigenvalues for parameter values outside of the training data. However, the diagonalization of the system employed to obtain this equation relies on the bi-orthogonality of the reduced-bases $\b{\Phi_m}(\mu)$ and $\b{\Psi_m}(\mu)$. We highlight that, for heat conduction problems, the direct eigenvectors and adjoint eigenvectors are identical because the operator $\b{A}(\mu)$ is self-adjoint. As a result, $\b{\Phi_m}(\mu)=\b{\Psi_m}(\mu)$, so adjoint eigenvectors can be disregarded. The bi-orthogonality of $\b{\Phi_m}(\mu)$ is satisfied for sampled parameter values $\mu_k$, but it's not exactly satisfied for unsampled parameters when eigenmodes are interpolated directly or by EDM-based interpolation. When constructing a ROM using this approach, our assumption is that bi-orthogonality is relatively well satisfied for the interpolated basis of eigenmodes, so that the system is diagonalized as in Eq.~\eqref{DiagonalizedSystem} where the dynamics operator is $\b{\Lambda_m}(\mu)$, a diagonal matrix of the interpolated eigenvalues.

Consequently, a trajectory from a given initial condition can be simulated for an unsampled parameter $\mu$, if $m$ eigenmodes in a reduced-order basis $\b{\Phi_m}(\mu) \in \mathbb{C}^{n \times m}$ and $m$ eigenvalues in the diagonal matrix $\b{\Lambda_m}(\mu) \in \mathbb{C}^{m \times m}$ are known at the parameter $\mu$. We use our reduced eigenmode representation based on EDMs, which requires the interpolation of $r$ EDM coefficients to compute the shape of eigenmodes for unsampled parameters. The proposed method is applied to build reduced representations of each one of the $m=6$ eigenmodes retaining $r=2$ EDMs in each case. For comparison, additional ROMs are built using a direct interpolation of the eigenmodes in physical space over the parameter interval, which involves interpolating $m$ vectors of dimension $n$. A cubic-spline interpolation is used for the eigenvalues to build the diagonal matrix $\b{\Lambda_m}(\mu)$. As a third alternative for comparison, we include the interpolation of solutions obtained by simulating multiple ROMs that are built for the training parameter values. This approach has been previously used to estimate a solution for an unsampled parameter using local ROMs previously constructed for sampled parameters~\cite{fonzi2020prsa}. Consequently, for the interpolation of ROM solutions, $p=4$ trajectories are simulated using the available ROMs at the sampled parameter values, and then their results are interpolated in physical space.

It is important to consider that the steady-state temperature field used as a linearization point for the battery problem depends on the parameterized boundary condition. Therefore, the definition of a perturbation from the equilibrium is also parameter-dependent, as shown in appendix~\ref{app:therm}. In order to build and simulate a ROM using Eq.~\eqref{ROMApproximation}, an equilibrium for the queried parameter, $\b{\bar x}(\mu)$, must be estimated or computed. For the results presented in this work, we assume that the equilibrium at the queried parameter can be pre-computed.  Importantly, for all three alternatives being compared, the same pre-computed equilibrium is used as a linearization point.

\begin{figure}
    \centering
    \includegraphics[width=17cm]{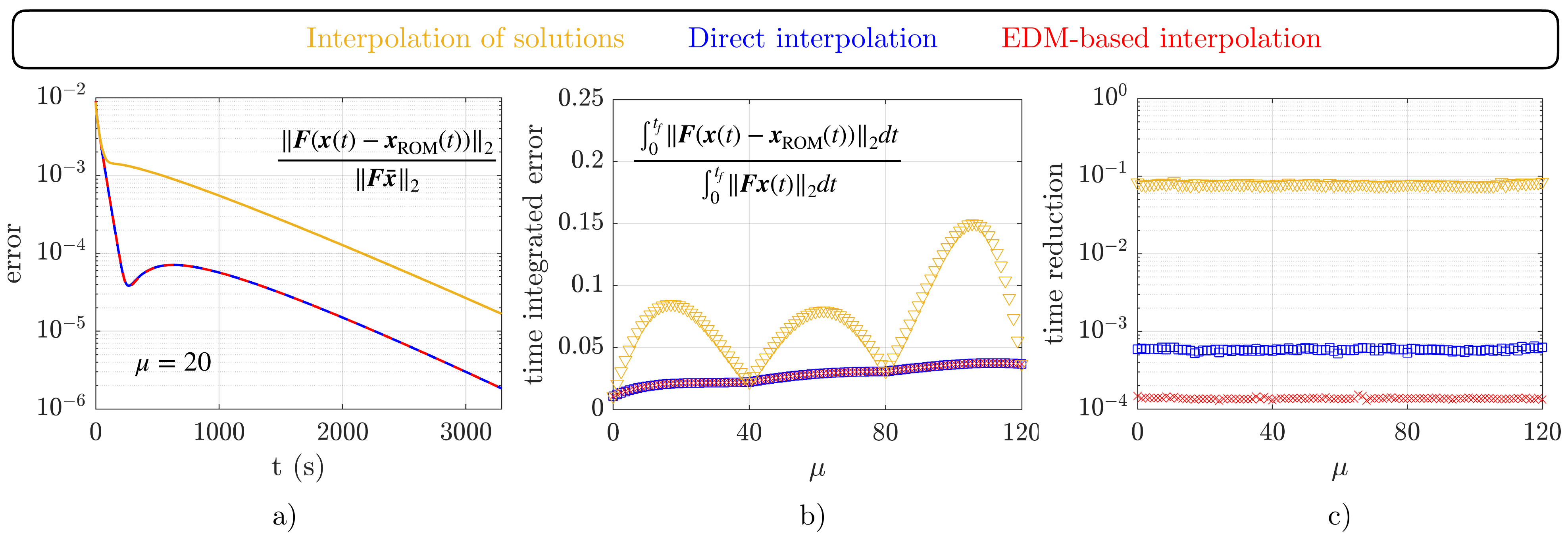}
    \caption{Performance of parameterized ROMs compared against the full-order model for the heat conduction in a battery problem. Projection-based ROMs using eigenmodes are built and three strategies to interpolate in parameter space from information at a few training parameter values are tested: (1) interpolating the solutions generated with available ROMs at the sampled parameter values, (2) directly interpolating eigenmodes in physical space, and (3) EDM-based interpolation, meaning that we interpolate the low-dimensional EDM coefficients. In the latter two cases, the interpolated eigenmodes are used for projection and interpolated eigenvalues are used to propagate the dynamics. (\textit{a}) Instantaneous error as a function of time for a single trajectory with $\mu=20$. (\textit{b}) Time integrated error over parameter space. (\textit{c}) Simulation time reduction.}
    \label{ROMResultsFigure}
\end{figure}

The performance of these three different approaches is compared using the error against the full-order model shown in Fig.~\ref{ROMResultsFigure} (\textit{a}). The initial condition is selected as the equilibrium for a parameter outside the range of sampled parameters. The results presented correspond to an initial condition set to the equilibrium at $\mu=300$, but similar results are obtained for tests conducted with parameters within the range of $\mu=120$ to $\mu=300$. Trajectories are simulated and 1000 uniformly spaced time steps are registered over a time horizon of 5 characteristic times based on the slowest eigenvalue for the first sampled parameter. For $\mu = 20$, which lies midway between two sampled parameter values, the ROMs built using directly interpolated eigenmodes and those using EDM-based interpolation outperform the interpolation of solutions strategy, as shown in Fig.~\ref{ROMResultsFigure} (\textit{a}). The performance of the ROMs is evaluated over a test set of parameter values that discretize the interval of interest between $\mu=0$ and $\mu=120$ using $100$ uniformly spaced steps. The normalized time-integrated error between the full-order model and the ROM trajectories as a function of the test parameter values is shown for the different interpolation strategies in Fig.~\ref{ROMResultsFigure} (\textit{b}). We find that direct and EDM-based interpolations share a similar performance, both outperforming the solution interpolation strategy at most parameter values, with the exception of the queried (training) parameter values. This is in agreement with the results presented in the previous subsections, with just two EDMs proving to be enough to capture the parameter-dependent deformation of the eigenmodes.


Lastly, Fig.~\ref{ROMResultsFigure} (\textit{c}) shows the reduction in simulation time, compared to the full order model, obtained with each of the ROMs based on the different interpolation strategies. The time is measured as CPU time in MATLAB, taking the average of 100 simulations per parameter to mitigate variations. For the solution interpolation approach, ROMs at each of the $p=4$  training parameter values need to be simulated, then the trajectories have to be interpolated in physical space for the testing parameter. For the other two approaches, eigenvalues and eigenmodes are interpolated and then a single ROM needs to be simulated for the testing parameter. For this example, we find that the ROM based on solution interpolation only reduces the computational time by a factor of 10 compared to the full order model, whereas the ROMs built using direct and EDM-based eigenmode interpolations reduce it by $3$ and $4$ orders of magnitude, respectively. Importantly, results obtained with the proposed approach leveraging EDMs show negligible loss of information compared to interpolating eigenmodes in physical space, while being an order of magnitude faster.

\section{\label{sec:conclusions}Conclusions}

In this work, we developed a data-driven framework to build interpretable low-order representations of eigenmodes in parameterized dynamical systems that capture eigenmode deformation through parameter space. The method takes as input a dataset of eigenmodes computed at a few queried parameter values across the interval of interest and relies on the singular value decomposition (SVD) to extract an orthogonal basis of \emph{eigen-deformation modes} (EDMs) that optimally captures the deformation of eigenmodes with the parameter. Furthermore, we show that, for a given eigenmode, EDMs are a discrete approximation of the eigenfunctions of the two-point spatial correlation tensor defined for an ensemble of instances of that eigenmode sampled across parameter space.

The method was tested on eigenmodes from linear or linearized problems arising in heat transfer, fluid mechanics, and structural vibrations, demonstrating its potential to effectively compress eigenmode parameter-dependence in a range of applications across civil, mechanical, and aerospace engineering. Importantly, we show that EDMs provide valuable insights into how parameter changes impact the eigenmodes and the behavior of the system under study. Furthermore, EDMs can be leveraged to significantly reduce the computational time required to build a projection-based reduced-order model for an unsampled parameter value. This is achieved by interpolating the reduced-order representations of each eigenmode at the queried parameter values instead of their full order counterparts in physical space. Moreover, the accuracy of the reduced-order representation can be controlled with the number of retained EDMs. 

The proposed method is data-driven and its implementation is straightforward, making it easy for practitioners to adopt, extend, and apply to a variety of fields where parameter-dependent eigenmode deformation is present. We remark that the compression capabilities of the method are presently limited by the use of linear dimensionality reduction. Future work could improve this by leveraging nonlinear dimensionality reduction techniques, however, achieving this without sacrificing the interpretability of EDMs is an open question. Another interesting avenue of future research is to apply the developed framework to capture the parameter dependence of modes arising from other — data-driven or equation-based —modal decompositions.

\begin{acknowledgments}
We gratefully acknowledge J. Lemus and D. Delgado for helpful comments and insightful discussions.
This work was funded by U. of Chile internal grant U-Inicia-003/21 and ANID Fondecyt 11220465.
\end{acknowledgments}


\providecommand{\noopsort}[1]{}\providecommand{\singleletter}[1]{#1}%

\appendix

\section{\label{app:therm}Thermal problems with internal heat generation}
Spatial discretization of transient heat conduction problems with internal heat generation can be expressed as a system of differential equations of the form

\begin{equation}
    \b{E} \b{\dot y}= \b{A}(\mu)\b{y} + \b{b},\label{eq:sys_source}
\end{equation}

\noindent where $\b{y}\in \mathbb{R}^{n}$ is the vector of temperatures at the $n$ mesh nodes and $\b{b}\in \mathbb{R}^{n}$ is the discretized internal heat generation term. $\b{E}\in \mathbb{R}^{n \times n}$ is the mass matrix and $\b{A}(\mu)\in \mathbb{R}^{n \times n}$ is the dynamics operator, which includes the boundary condition parameter dependence.

Let $\b{\bar y}$ be an equilibrium of the system for a given parameter value $\mu$, such that

\begin{equation}
    \b{0} = \b{A}(\mu) \b{\bar y} + \b{b}.\label{eq:ybar}
\end{equation}
Then, if we solve for $\b{b}$ in Eq.~\eqref{eq:ybar} and substitute it into Eq.~\eqref{eq:sys_source}, we can recognize that the dynamics of a perturbation away from the equilibrium $\b{x}=\b{y} - \b{\bar y}$ are governed by

\begin{equation}
\b{E} \b{\dot x} = \b{A}(\mu) \b{x},
\end{equation}
that has the exact same form of Eq.~\eqref{dynamicalSystem}, so eigenmodes and eigenvalues can be obtained by solving the corresponding eigenproblem for different parameter values.

\section{\label{app:mech}Mechanical systems in first order form}

Mechanical systems governed by linear elasticity equations in the absence of external forces can be expressed, when the spatial domain is discretized, as a system of second order differential equations

\begin{equation}
    \b{M \ddot y} = \b{K}(\mu) \b{y},\label{eq:sys_mech}
\end{equation}
\noindent where $\b{y}\in \mathbb{R}^{n}$ is the vector of displacements, the overdot denotes time differentiation, and $\b{M}\in \mathbb{R}^{n \times n}$ and $\b{K}(\mu)\in \mathbb{R}^{n\times n}$ are the mass and stiffness matrices respectively, with the latter incorporating the parameter dependence.

Let $\b{x}\in \mathbb{R}^{2 n}$ be an augmented state variable obtained by concatenating the displacements and velocities as follows

\begin{equation}
    \b{x} = \begin{bmatrix}
           \b{y} \\
           \b{\dot y}
         \end{bmatrix}.
\end{equation}
We may now rewrite Eq.~\eqref{eq:sys_mech} as a system of first order differential equations of the form
\begin{equation}\label{mechFirstOrder}
    \begin{bmatrix}
           \b{I} & \b{0}\\
           \b{0} & \b{M}
         \end{bmatrix}
          \begin{bmatrix}
           \b{\dot y} \\
           \b{\ddot y}
         \end{bmatrix} =
        \begin{bmatrix}
           \b{0} & \b{I}\\
           \b{K} & \b{0}
         \end{bmatrix}
          \begin{bmatrix}
           \b{y} \\
           \b{\dot y}
         \end{bmatrix},
\end{equation}

\noindent where $\b{I} \in \mathbb{R}^{n\times n}$ and $\b{0} \in \mathbb{R}^{n\times n}$ are the identity matrix, and a matrix of zeros, respectively. The above is a linear system in the form of Eq.~\eqref{dynamicalSystem}, for the augmented state variable $\b{x}$, where

\begin{equation}\label{mechMatrices}
    \b{E} = \begin{bmatrix}
           \b{I} & \b{0}\\
           \b{0} & \b{M}
         \end{bmatrix}\:\: \textrm{and} \:\: \b{A} = 
         \begin{bmatrix}
           \b{0} & \b{I}\\
           \b{K} & \b{0}
         \end{bmatrix},
\end{equation}

\noindent so the corresponding generalized eigenproblem with the matrices in~\eqref{mechMatrices} can be solved to obtain eigenmodes and eigenvalues.

\section{\label{app:preProcessing}Consistent normalization}
 
 A consistent normalization of eigenmodes is critical for the computation of EDMs with the proposed method. The usual convention is to normalize eigenmodes so that they have unit norm based on the inner-product considering the mass matrix $\b{E}$, that is $\b{\phi_i^{T}}(\mu_k)\b{E}\b{\phi_i}(\mu_k)=1$. However, even after normalization, the phase of the complex-valued eigenmodes (or their sign if they are real-valued) is completely arbitrary and needs to be aligned between mode instances obtained for different parameter values. In this work, we address this with a simple pre-processing of the data to maintain a consistent normalization and phase (or sign) alignment across the sampled parameter values. This ensures the smoothness of the EDM coefficients within the parameter interval of interest, as shown in section~\ref{sec:results}. The process used to align real-valued modes is shown in algorithm~\ref{algo_real}.


\begin{algorithm}[H]
    \caption{Pre-processing of real-valued modes} 
    \begin{flushleft}
    \textbf{Inputs:} $\b{\phi}_i^{(k)}$ for $i = 1, \ldots, m$ and $k = 1, \ldots, p$ \\
    \textbf{Outputs:} Sign-aligned modes $\tilde{\b{\phi}}_i^{(k)}$ \\
    \end{flushleft}
    \begin{algorithmic}[1]
        \For {$i=1,2,\ldots, m$}
            \State $\tilde{\b{\phi}}_i^{(1)} \leftarrow \b{\phi}_i^{(1)}$
            \For {$k=1,2,\ldots, p-1$}
                \State $a_{ik}\leftarrow (\b{\phi}_i^{(k)} )^{T}\b{E}\b{\phi}_i^{(k+1)}$
                \If{$a_{ik} \leq 0$} 
                    \State $\tilde{\b{\phi}}_i^{(k+1)} \leftarrow - \b{\phi}_i^{(k+1)}$
                \Else
                    \State $\tilde{\b{\phi}}_i^{(k+1)} \leftarrow   \b{\phi}_i^{(k+1)}$
                \EndIf
            \EndFor
        \EndFor
    \end{algorithmic} \label{algo_real}
\end{algorithm}

For complex-valued eigenmodes, an optimization problem was formulated to find the angle $\theta_{ik}$ that best aligns the phase of the $i$-th eigenmode for parameter $\mu_k$, $\b{\phi}_i^{(k)}$, with the same eigenmode but sampled for the first parameter value $\b{\phi}_i^{(1)}$
\begin{equation}\label{alignOpt}
    \theta_{ik} = \text{arg}\,\min\limits_{\theta}\, \|\b{F}( \b{\phi}_i^k e^{j \theta} - \b{\phi}_i^1)\|_2,
\end{equation}
where we use $j$ for the imaginary unit here. The resulting process used to align complex-valued eigenmodes by iteratively solving this optimization problem is shown in algorithm~\ref{algo_complex}.

\begin{algorithm}[H]
    \caption{Pre-processing of complex-valued modes}
    \begin{flushleft}
    \textbf{Inputs:} $\b{\phi}_i^{(k)}$ for $i = 1, \ldots, m$ and $k = 1, \ldots, p$ \\
    \textbf{Outputs:} Phase-aligned modes $\tilde{\b{\phi}}_i^{(k)}$ \\
    \end{flushleft}
    \begin{algorithmic}[1]
        \For {$i=1,2,\ldots, m$}
            \State $\tilde{\b{\phi}}_i^{(1)} \leftarrow \b{\phi}_i^{(1)}$
            \For {$k=2,3,\ldots, p$}
                \State $\theta_{ik} \leftarrow \text{arg}\,\min\limits_{\theta}\, \| \b{F}(\b{\phi}_i^{(k)} e^{j \theta} - \b{\phi}_i^{(1)})\|_2$
                \State $\tilde{\b{\phi}}_i^{(k)} \leftarrow \b{\phi}_i^{(k)} \theta_{ik}$
            \EndFor
        \EndFor
    \end{algorithmic}\label{algo_complex}
\end{algorithm}

The corresponding of the above pre-processing steps was applied to the eigenmode datasets before computing EDMs in the examples presented to ensure a consistent normalization and phase or sign alignment.

\section{Battery modeling parameters}\label{app:BatteryParameters}
\begin{figure}[h]
    \centering
    \includegraphics[width=8cm]{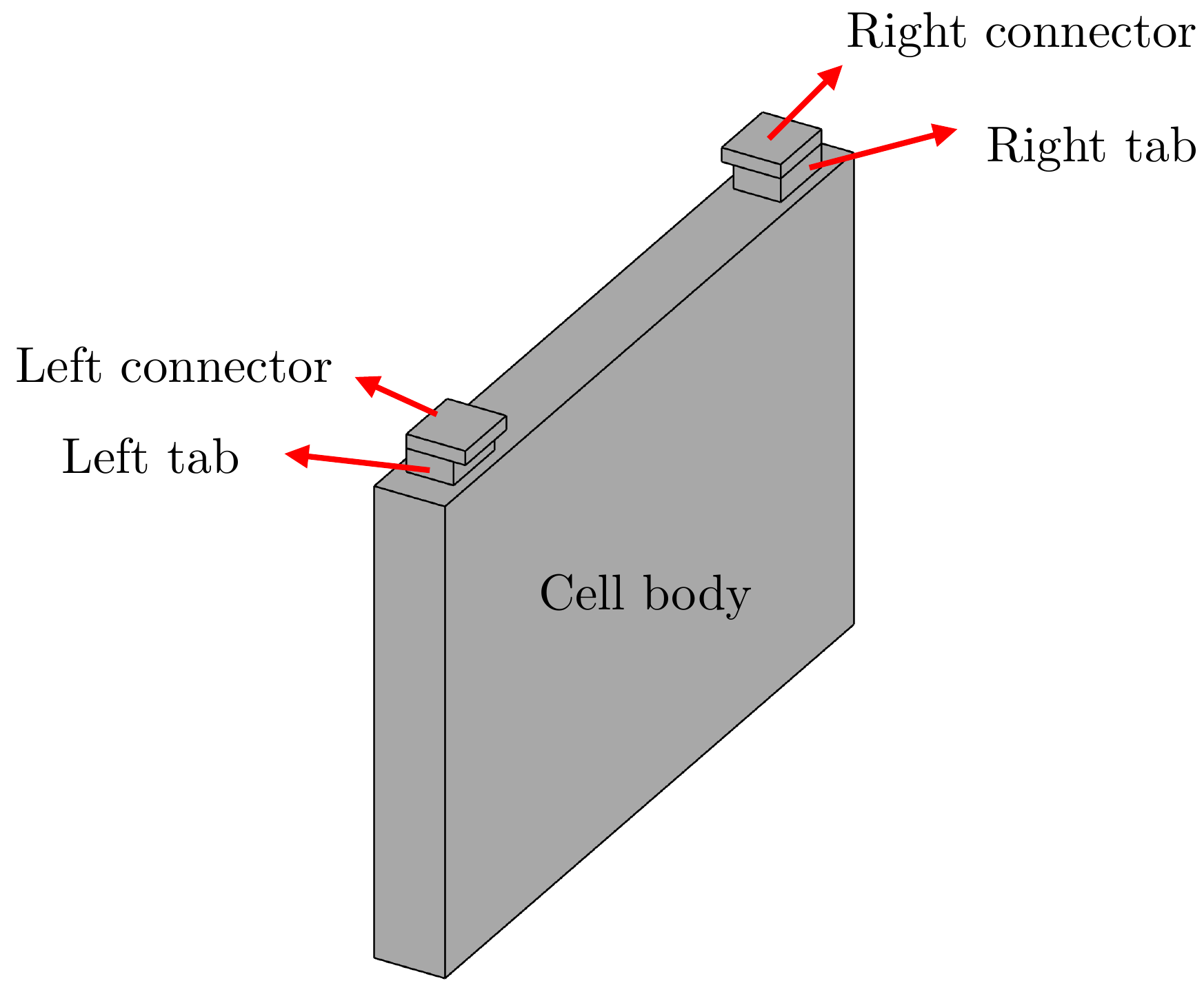}
    \caption{Cells created using MATLAB's PDE toolbox. The five cell volumes are shown: two tabs, two connectors and the cell body.}
    \label{batteryCellGeometry}
\end{figure}

The creation of the modular prismatic battery geometry used in the heat transfer example was done using \textit{createBatteryModuleGeometry} function in MATLAB's PDE toolbox. This function creates a geometry consisting of a specified number of cells. Each cell is divided into 5 volumes: the cell body, 2 connectors and 2 tabs, as shown in Fig~\ref{batteryCellGeometry}. It is possible to assign different dimensions and material properties to each one of the volumes. Table~\ref{batteryParametersTable} shows all the fixed parameters used to model the battery.

\begin{table}[h]\label{batteryParametersTable}
\centering
\begin{tabular}[t]{lc@{\hskip 0.5in}c}
\hline
&Value & Unit\\
\hline
Cell width                                        &0,015           & $\si{\m}$\\
Cell thickness                                    &0,015           & $\si{\m}$\\
Tab thickness                                     &0,01            & $\si{\m}$\\
Tab width                                         &0,015           & $\si{\m}$\\
Cell height                                       &0,1             & $\si{\m}$\\
Tab height                                        &0,005           & $\si{\m}$\\
Connector height                                  &0,003           & $\si{\m}$\\
Number of cells                                   &20              & -\\
In plane cell thermal conductivity                &80              & $\si{\W / (\K \m)}$\\
Through plane cell thermal   conductivity         &2               & $\si{\W / (\K \m)}$\\
Tab thermal conductivity                          &386             & $\si{\W / (\K \m)}$\\
Connector thermal conductivity                    &400             & $\si{\W / (\K \m)}$\\
Cell mass density                                 &780             & $\si{\kilogram / \m ^3}$\\
Tabs mass density                                 &2700            & $\si{\kilogram / \m ^3}$\\
Connectors mass density                           &540             & $\si{\kilogram / \m ^3}$\\
Cell specific heat                                &785             & $\si{\J / (\kilogram \K)}$\\
Tabs specific heat                                &890             & $\si{\J / (\kilogram \K)}$\\
Connectors specific heat                          &840             & $\si{\J / (\kilogram \K)}$\\
Ambient temperature  $T_{\infty}$                 &293             & $\si{\K}$\\
Heat generation per cell                          &15              & $\si{\W}$\\
\hline
\end{tabular}
\caption{Fixed parameters used for modular battery modeling.}
\end{table}%

\section{Third eigenmode interpolation}\label{app:ThirdEigenmode}
\begin{figure}
    \centering
    \includegraphics[width=17cm]{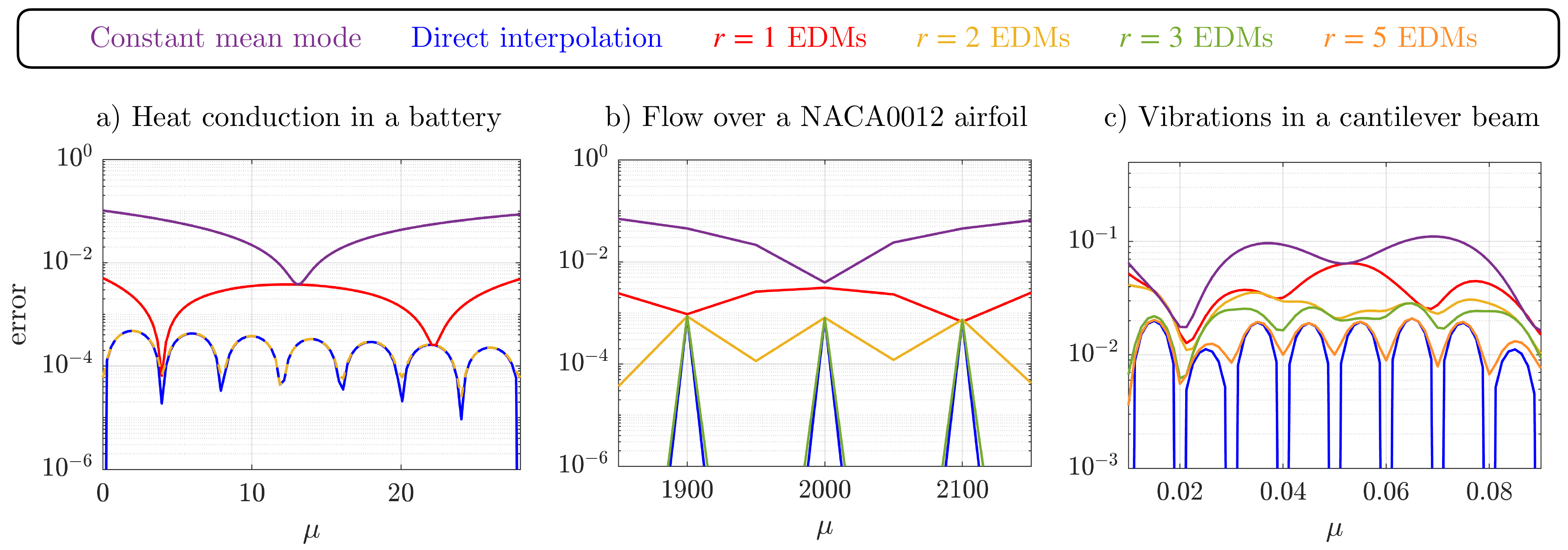}
    \caption{Spatially integrated error between the third eigenmode and its interpolation over parameter space from a few sampled parameter values for the (\textit{a}) heat transfer, (\textit{b}) fluid mechanics, and (\textit{c}) structural dynamics examples. Direct interpolation of the eigenmodes in physical space is compared to EDM-based interpolation using a different number of EDMs.}
    \label{3rd_eigenmode}
\end{figure}

To complement the results discussed in section~\S\ref{sec:Interpolation}, here we present the direct and EDM-based interpolation of an additional eigenmode (beyond the leading one) for each of our examples. Importantly, the same trends observed for the first eigenmode are found here for the third eigenmode of each example, as shown in Fig.~\ref{3rd_eigenmode}. This provides further evidence that EDMs can effectively compress eigenmode deformation across parameter space and produce an interpretable low-order representation that is as accurate as the interpolation in physical space.

\begin{thebibliography}{51}%
\makeatletter
\providecommand \@ifxundefined [1]{%
 \@ifx{#1\undefined}
}%
\providecommand \@ifnum [1]{%
 \ifnum #1\expandafter \@firstoftwo
 \else \expandafter \@secondoftwo
 \fi
}%
\providecommand \@ifx [1]{%
 \ifx #1\expandafter \@firstoftwo
 \else \expandafter \@secondoftwo
 \fi
}%
\providecommand \natexlab [1]{#1}%
\providecommand \enquote  [1]{``#1''}%
\providecommand \bibnamefont  [1]{#1}%
\providecommand \bibfnamefont [1]{#1}%
\providecommand \citenamefont [1]{#1}%
\providecommand \href@noop [0]{\@secondoftwo}%
\providecommand \href [0]{\begingroup \@sanitize@url \@href}%
\providecommand \@href[1]{\@@startlink{#1}\@@href}%
\providecommand \@@href[1]{\endgroup#1\@@endlink}%
\providecommand \@sanitize@url [0]{\catcode `\\12\catcode `\$12\catcode `\&12\catcode `\#12\catcode `\^12\catcode `\_12\catcode `\%12\relax}%
\providecommand \@@startlink[1]{}%
\providecommand \@@endlink[0]{}%
\providecommand \url  [0]{\begingroup\@sanitize@url \@url }%
\providecommand \@url [1]{\endgroup\@href {#1}{\urlprefix }}%
\providecommand \urlprefix  [0]{URL }%
\providecommand \Eprint [0]{\href }%
\providecommand \doibase [0]{https://doi.org/}%
\providecommand \selectlanguage [0]{\@gobble}%
\providecommand \bibinfo  [0]{\@secondoftwo}%
\providecommand \bibfield  [0]{\@secondoftwo}%
\providecommand \translation [1]{[#1]}%
\providecommand \BibitemOpen [0]{}%
\providecommand \bibitemStop [0]{}%
\providecommand \bibitemNoStop [0]{.\EOS\space}%
\providecommand \EOS [0]{\spacefactor3000\relax}%
\providecommand \BibitemShut  [1]{\csname bibitem#1\endcsname}%
\let\auto@bib@innerbib\@empty
\bibitem [{\citenamefont {Brunton}\ and\ \citenamefont {Noack}(2015)}]{brunton2015amr}%
  \BibitemOpen
  \bibfield  {author} {\bibinfo {author} {\bibfnamefont {S.~L.}\ \bibnamefont {Brunton}}\ and\ \bibinfo {author} {\bibfnamefont {B.~R.}\ \bibnamefont {Noack}},\ }\bibfield  {title} {\bibinfo {title} {Closed-loop turbulence control: Progress and challenges},\ }\href@noop {} {\bibfield  {journal} {\bibinfo  {journal} {Applied Mechanics Reviews}\ }\textbf {\bibinfo {volume} {67}},\ \bibinfo {pages} {050801} (\bibinfo {year} {2015})}\BibitemShut {NoStop}%
\bibitem [{\citenamefont {Herrmann}\ \emph {et~al.}(2020)\citenamefont {Herrmann}, \citenamefont {Behzad}, \citenamefont {Cardemil}, \citenamefont {Calder{\'o}n-Mu{\~n}oz},\ and\ \citenamefont {Fern{\'a}ndez}}]{herrmann2020se}%
  \BibitemOpen
  \bibfield  {author} {\bibinfo {author} {\bibfnamefont {B.}~\bibnamefont {Herrmann}}, \bibinfo {author} {\bibfnamefont {M.}~\bibnamefont {Behzad}}, \bibinfo {author} {\bibfnamefont {J.~M.}\ \bibnamefont {Cardemil}}, \bibinfo {author} {\bibfnamefont {W.~R.}\ \bibnamefont {Calder{\'o}n-Mu{\~n}oz}},\ and\ \bibinfo {author} {\bibfnamefont {R.~M.}\ \bibnamefont {Fern{\'a}ndez}},\ }\bibfield  {title} {\bibinfo {title} {Conjugate heat transfer model for feedback control and state estimation in a volumetric solar receiver},\ }\href@noop {} {\bibfield  {journal} {\bibinfo  {journal} {Solar Energy}\ }\textbf {\bibinfo {volume} {198}},\ \bibinfo {pages} {343} (\bibinfo {year} {2020})}\BibitemShut {NoStop}%
\bibitem [{\citenamefont {Herrmann}\ \emph {et~al.}(2022)\citenamefont {Herrmann}, \citenamefont {Brunton}, \citenamefont {Pohl},\ and\ \citenamefont {Semaan}}]{herrmann2022prf}%
  \BibitemOpen
  \bibfield  {author} {\bibinfo {author} {\bibfnamefont {B.}~\bibnamefont {Herrmann}}, \bibinfo {author} {\bibfnamefont {S.~L.}\ \bibnamefont {Brunton}}, \bibinfo {author} {\bibfnamefont {J.~E.}\ \bibnamefont {Pohl}},\ and\ \bibinfo {author} {\bibfnamefont {R.}~\bibnamefont {Semaan}},\ }\bibfield  {title} {\bibinfo {title} {Gust mitigation through closed-loop control. {II}. {F}eedforward and feedback control},\ }\href@noop {} {\bibfield  {journal} {\bibinfo  {journal} {Physical Review Fluids}\ }\textbf {\bibinfo {volume} {7}},\ \bibinfo {pages} {024706} (\bibinfo {year} {2022})}\BibitemShut {NoStop}%
\bibitem [{\citenamefont {Gunzburger}(1999)}]{gunzburger1999ijnmf}%
  \BibitemOpen
  \bibfield  {author} {\bibinfo {author} {\bibfnamefont {M.~D.}\ \bibnamefont {Gunzburger}},\ }\bibfield  {title} {\bibinfo {title} {Sensitivities, adjoints and flow optimization},\ }\href@noop {} {\bibfield  {journal} {\bibinfo  {journal} {International Journal for Numerical Methods in Fluids}\ }\textbf {\bibinfo {volume} {31}},\ \bibinfo {pages} {53} (\bibinfo {year} {1999})}\BibitemShut {NoStop}%
\bibitem [{\citenamefont {Giles}\ and\ \citenamefont {Pierce}(2000)}]{giles2000ftc}%
  \BibitemOpen
  \bibfield  {author} {\bibinfo {author} {\bibfnamefont {M.~B.}\ \bibnamefont {Giles}}\ and\ \bibinfo {author} {\bibfnamefont {N.~A.}\ \bibnamefont {Pierce}},\ }\bibfield  {title} {\bibinfo {title} {An introduction to the adjoint approach to design},\ }\href@noop {} {\bibfield  {journal} {\bibinfo  {journal} {Flow, Turbulence and Combustion}\ }\textbf {\bibinfo {volume} {65}},\ \bibinfo {pages} {393} (\bibinfo {year} {2000})}\BibitemShut {NoStop}%
\bibitem [{\citenamefont {Herrmann-Priesnitz}\ \emph {et~al.}(2016)\citenamefont {Herrmann-Priesnitz}, \citenamefont {Calder{\'o}n-Mu{\~n}oz}, \citenamefont {Valencia},\ and\ \citenamefont {Soto}}]{herrmann2016ate}%
  \BibitemOpen
  \bibfield  {author} {\bibinfo {author} {\bibfnamefont {B.}~\bibnamefont {Herrmann-Priesnitz}}, \bibinfo {author} {\bibfnamefont {W.~R.}\ \bibnamefont {Calder{\'o}n-Mu{\~n}oz}}, \bibinfo {author} {\bibfnamefont {A.}~\bibnamefont {Valencia}},\ and\ \bibinfo {author} {\bibfnamefont {R.}~\bibnamefont {Soto}},\ }\bibfield  {title} {\bibinfo {title} {Thermal design exploration of a swirl flow microchannel heat sink for high heat flux applications based on numerical simulations},\ }\href@noop {} {\bibfield  {journal} {\bibinfo  {journal} {Applied Thermal Engineering}\ }\textbf {\bibinfo {volume} {109}},\ \bibinfo {pages} {22} (\bibinfo {year} {2016})}\BibitemShut {NoStop}%
\bibitem [{\citenamefont {Herrmann-Priesnitz}\ \emph {et~al.}(2018)\citenamefont {Herrmann-Priesnitz}, \citenamefont {Calder{\'o}n-Mu{\~n}oz}, \citenamefont {Diaz},\ and\ \citenamefont {Soto}}]{herrmann2018hmt}%
  \BibitemOpen
  \bibfield  {author} {\bibinfo {author} {\bibfnamefont {B.}~\bibnamefont {Herrmann-Priesnitz}}, \bibinfo {author} {\bibfnamefont {W.~R.}\ \bibnamefont {Calder{\'o}n-Mu{\~n}oz}}, \bibinfo {author} {\bibfnamefont {G.}~\bibnamefont {Diaz}},\ and\ \bibinfo {author} {\bibfnamefont {R.}~\bibnamefont {Soto}},\ }\bibfield  {title} {\bibinfo {title} {Heat transfer enhancement strategies in a swirl flow minichannel heat sink based on hydrodynamic receptivity},\ }\href@noop {} {\bibfield  {journal} {\bibinfo  {journal} {International Journal of Heat and Mass Transfer}\ }\textbf {\bibinfo {volume} {127}},\ \bibinfo {pages} {245} (\bibinfo {year} {2018})}\BibitemShut {NoStop}%
\bibitem [{\citenamefont {Grieves}\ and\ \citenamefont {Vickers}(2017)}]{grieves2017springer}%
  \BibitemOpen
  \bibfield  {author} {\bibinfo {author} {\bibfnamefont {M.}~\bibnamefont {Grieves}}\ and\ \bibinfo {author} {\bibfnamefont {J.}~\bibnamefont {Vickers}},\ }\bibfield  {title} {\bibinfo {title} {Digital twin: Mitigating unpredictable, undesirable emergent behavior in complex systems},\ }in\ \href@noop {} {\emph {\bibinfo {booktitle} {Transdisciplinary perspectives on complex systems: New findings and approaches}}}\ (\bibinfo  {publisher} {Springer},\ \bibinfo {year} {2017})\ pp.\ \bibinfo {pages} {85--113}\BibitemShut {NoStop}%
\bibitem [{\citenamefont {Hartmann}\ \emph {et~al.}(2018)\citenamefont {Hartmann}, \citenamefont {Herz},\ and\ \citenamefont {Wever}}]{hartmann2018springer}%
  \BibitemOpen
  \bibfield  {author} {\bibinfo {author} {\bibfnamefont {D.}~\bibnamefont {Hartmann}}, \bibinfo {author} {\bibfnamefont {M.}~\bibnamefont {Herz}},\ and\ \bibinfo {author} {\bibfnamefont {U.}~\bibnamefont {Wever}},\ }\bibfield  {title} {\bibinfo {title} {Model order reduction a key technology for digital twins},\ }in\ \href@noop {} {\emph {\bibinfo {booktitle} {Reduced-order modeling (ROM) for simulation and optimization}}}\ (\bibinfo  {publisher} {Springer},\ \bibinfo {year} {2018})\ pp.\ \bibinfo {pages} {167--179}\BibitemShut {NoStop}%
\bibitem [{\citenamefont {Niederer}\ \emph {et~al.}(2021)\citenamefont {Niederer}, \citenamefont {Sacks}, \citenamefont {Girolami},\ and\ \citenamefont {Willcox}}]{niederer2021natcs}%
  \BibitemOpen
  \bibfield  {author} {\bibinfo {author} {\bibfnamefont {S.~A.}\ \bibnamefont {Niederer}}, \bibinfo {author} {\bibfnamefont {M.~S.}\ \bibnamefont {Sacks}}, \bibinfo {author} {\bibfnamefont {M.}~\bibnamefont {Girolami}},\ and\ \bibinfo {author} {\bibfnamefont {K.}~\bibnamefont {Willcox}},\ }\bibfield  {title} {\bibinfo {title} {Scaling digital twins from the artisanal to the industrial},\ }\href@noop {} {\bibfield  {journal} {\bibinfo  {journal} {Nature Computational Science}\ }\textbf {\bibinfo {volume} {1}},\ \bibinfo {pages} {313} (\bibinfo {year} {2021})}\BibitemShut {NoStop}%
\bibitem [{\citenamefont {Fu}\ and\ \citenamefont {He}(2001)}]{fu2001book}%
  \BibitemOpen
  \bibfield  {author} {\bibinfo {author} {\bibfnamefont {Z.-F.}\ \bibnamefont {Fu}}\ and\ \bibinfo {author} {\bibfnamefont {J.}~\bibnamefont {He}},\ }\href@noop {} {\emph {\bibinfo {title} {Modal analysis}}}\ (\bibinfo  {publisher} {Elsevier},\ \bibinfo {year} {2001})\BibitemShut {NoStop}%
\bibitem [{\citenamefont {Bergman}\ \emph {et~al.}(2011)\citenamefont {Bergman}, \citenamefont {Lavine}, \citenamefont {Incropera},\ and\ \citenamefont {DeWitt}}]{bergman2011book}%
  \BibitemOpen
  \bibfield  {author} {\bibinfo {author} {\bibfnamefont {T.~L.}\ \bibnamefont {Bergman}}, \bibinfo {author} {\bibfnamefont {A.~S.}\ \bibnamefont {Lavine}}, \bibinfo {author} {\bibfnamefont {F.~P.}\ \bibnamefont {Incropera}},\ and\ \bibinfo {author} {\bibfnamefont {D.~P.}\ \bibnamefont {DeWitt}},\ }\href@noop {} {\emph {\bibinfo {title} {Introduction to heat transfer}}}\ (\bibinfo  {publisher} {John Wiley \& Sons},\ \bibinfo {year} {2011})\BibitemShut {NoStop}%
\bibitem [{\citenamefont {Theofilis}(2011)}]{theofilis2011arfm}%
  \BibitemOpen
  \bibfield  {author} {\bibinfo {author} {\bibfnamefont {V.}~\bibnamefont {Theofilis}},\ }\bibfield  {title} {\bibinfo {title} {Global linear instability},\ }\href@noop {} {\bibfield  {journal} {\bibinfo  {journal} {Annual Review of Fluid Mechanics}\ }\textbf {\bibinfo {volume} {43}},\ \bibinfo {pages} {319} (\bibinfo {year} {2011})}\BibitemShut {NoStop}%
\bibitem [{\citenamefont {Taira}\ \emph {et~al.}(2017)\citenamefont {Taira}, \citenamefont {Brunton}, \citenamefont {Dawson}, \citenamefont {Rowley}, \citenamefont {Colonius}, \citenamefont {McKeon}, \citenamefont {Schmidt}, \citenamefont {Gordeyev}, \citenamefont {Theofilis},\ and\ \citenamefont {Ukeiley}}]{taira2017aiaa}%
  \BibitemOpen
  \bibfield  {author} {\bibinfo {author} {\bibfnamefont {K.}~\bibnamefont {Taira}}, \bibinfo {author} {\bibfnamefont {S.~L.}\ \bibnamefont {Brunton}}, \bibinfo {author} {\bibfnamefont {S.}~\bibnamefont {Dawson}}, \bibinfo {author} {\bibfnamefont {C.~W.}\ \bibnamefont {Rowley}}, \bibinfo {author} {\bibfnamefont {T.}~\bibnamefont {Colonius}}, \bibinfo {author} {\bibfnamefont {B.~J.}\ \bibnamefont {McKeon}}, \bibinfo {author} {\bibfnamefont {O.~T.}\ \bibnamefont {Schmidt}}, \bibinfo {author} {\bibfnamefont {S.}~\bibnamefont {Gordeyev}}, \bibinfo {author} {\bibfnamefont {V.}~\bibnamefont {Theofilis}},\ and\ \bibinfo {author} {\bibfnamefont {L.~S.}\ \bibnamefont {Ukeiley}},\ }\bibfield  {title} {\bibinfo {title} {Modal analysis of fluid flows: An overview},\ }\href@noop {} {\bibfield  {journal} {\bibinfo  {journal} {AIAA Journal}\ }\textbf {\bibinfo {volume} {55}},\ \bibinfo {pages} {4013} (\bibinfo {year} {2017})}\BibitemShut {NoStop}%
\bibitem [{\citenamefont {Taira}\ \emph {et~al.}(2020)\citenamefont {Taira}, \citenamefont {Hemati}, \citenamefont {Brunton}, \citenamefont {Sun}, \citenamefont {Duraisamy}, \citenamefont {Bagheri}, \citenamefont {Dawson},\ and\ \citenamefont {Yeh}}]{taira2019aiaa}%
  \BibitemOpen
  \bibfield  {author} {\bibinfo {author} {\bibfnamefont {K.}~\bibnamefont {Taira}}, \bibinfo {author} {\bibfnamefont {M.~S.}\ \bibnamefont {Hemati}}, \bibinfo {author} {\bibfnamefont {S.~L.}\ \bibnamefont {Brunton}}, \bibinfo {author} {\bibfnamefont {Y.}~\bibnamefont {Sun}}, \bibinfo {author} {\bibfnamefont {K.}~\bibnamefont {Duraisamy}}, \bibinfo {author} {\bibfnamefont {S.}~\bibnamefont {Bagheri}}, \bibinfo {author} {\bibfnamefont {S.~T.~M.}\ \bibnamefont {Dawson}},\ and\ \bibinfo {author} {\bibfnamefont {C.-A.}\ \bibnamefont {Yeh}},\ }\bibfield  {title} {\bibinfo {title} {Modal analysis of fluid flows: Applications and outlook},\ }\href@noop {} {\bibfield  {journal} {\bibinfo  {journal} {AIAA Journal}\ }\textbf {\bibinfo {volume} {58}},\ \bibinfo {pages} {998} (\bibinfo {year} {2020})}\BibitemShut {NoStop}%
\bibitem [{\citenamefont {Herrmann}\ \emph {et~al.}(2021)\citenamefont {Herrmann}, \citenamefont {Baddoo}, \citenamefont {Semaan}, \citenamefont {Brunton},\ and\ \citenamefont {McKeon}}]{herrmann2021jfm}%
  \BibitemOpen
  \bibfield  {author} {\bibinfo {author} {\bibfnamefont {B.}~\bibnamefont {Herrmann}}, \bibinfo {author} {\bibfnamefont {P.~J.}\ \bibnamefont {Baddoo}}, \bibinfo {author} {\bibfnamefont {R.}~\bibnamefont {Semaan}}, \bibinfo {author} {\bibfnamefont {S.~L.}\ \bibnamefont {Brunton}},\ and\ \bibinfo {author} {\bibfnamefont {B.~J.}\ \bibnamefont {McKeon}},\ }\bibfield  {title} {\bibinfo {title} {Data-driven resolvent analysis},\ }\href@noop {} {\bibfield  {journal} {\bibinfo  {journal} {Journal of Fluid Mechanics}\ }\textbf {\bibinfo {volume} {918}},\ \bibinfo {pages} {A10} (\bibinfo {year} {2021})}\BibitemShut {NoStop}%
\bibitem [{\citenamefont {Schmid}(2022)}]{schmid2022arfm}%
  \BibitemOpen
  \bibfield  {author} {\bibinfo {author} {\bibfnamefont {P.~J.}\ \bibnamefont {Schmid}},\ }\bibfield  {title} {\bibinfo {title} {Dynamic mode decomposition and its variants},\ }\href@noop {} {\bibfield  {journal} {\bibinfo  {journal} {Annual Review of Fluid Mechanics}\ }\textbf {\bibinfo {volume} {54}},\ \bibinfo {pages} {225} (\bibinfo {year} {2022})}\BibitemShut {NoStop}%
\bibitem [{\citenamefont {Baddoo}\ \emph {et~al.}(2023)\citenamefont {Baddoo}, \citenamefont {Herrmann}, \citenamefont {McKeon}, \citenamefont {Nathan~Kutz},\ and\ \citenamefont {Brunton}}]{baddoo2023prsa}%
  \BibitemOpen
  \bibfield  {author} {\bibinfo {author} {\bibfnamefont {P.~J.}\ \bibnamefont {Baddoo}}, \bibinfo {author} {\bibfnamefont {B.}~\bibnamefont {Herrmann}}, \bibinfo {author} {\bibfnamefont {B.~J.}\ \bibnamefont {McKeon}}, \bibinfo {author} {\bibfnamefont {J.}~\bibnamefont {Nathan~Kutz}},\ and\ \bibinfo {author} {\bibfnamefont {S.~L.}\ \bibnamefont {Brunton}},\ }\bibfield  {title} {\bibinfo {title} {Physics-informed dynamic mode decomposition},\ }\href@noop {} {\bibfield  {journal} {\bibinfo  {journal} {Proceedings of the Royal Society A}\ }\textbf {\bibinfo {volume} {479}},\ \bibinfo {pages} {20220576} (\bibinfo {year} {2023})}\BibitemShut {NoStop}%
\bibitem [{\citenamefont {Benner}\ \emph {et~al.}(2015)\citenamefont {Benner}, \citenamefont {Gugercin},\ and\ \citenamefont {Willcox}}]{benner2015sr}%
  \BibitemOpen
  \bibfield  {author} {\bibinfo {author} {\bibfnamefont {P.}~\bibnamefont {Benner}}, \bibinfo {author} {\bibfnamefont {S.}~\bibnamefont {Gugercin}},\ and\ \bibinfo {author} {\bibfnamefont {K.}~\bibnamefont {Willcox}},\ }\bibfield  {title} {\bibinfo {title} {A survey of projection-based model reduction methods for parametric dynamical systems},\ }\href@noop {} {\bibfield  {journal} {\bibinfo  {journal} {SIAM Review}\ }\textbf {\bibinfo {volume} {57}},\ \bibinfo {pages} {483} (\bibinfo {year} {2015})}\BibitemShut {NoStop}%
\bibitem [{\citenamefont {Rowley}\ and\ \citenamefont {Dawson}(2017)}]{rowley2017arfm}%
  \BibitemOpen
  \bibfield  {author} {\bibinfo {author} {\bibfnamefont {C.~W.}\ \bibnamefont {Rowley}}\ and\ \bibinfo {author} {\bibfnamefont {S.~T.}\ \bibnamefont {Dawson}},\ }\bibfield  {title} {\bibinfo {title} {Model reduction for flow analysis and control},\ }\href@noop {} {\bibfield  {journal} {\bibinfo  {journal} {Annual Review of Fluid Mechanics}\ }\textbf {\bibinfo {volume} {49}},\ \bibinfo {pages} {387} (\bibinfo {year} {2017})}\BibitemShut {NoStop}%
\bibitem [{\citenamefont {Antoulas}(2005)}]{antoulasbook}%
  \BibitemOpen
  \bibfield  {author} {\bibinfo {author} {\bibfnamefont {A.~C.}\ \bibnamefont {Antoulas}},\ }\href@noop {} {\emph {\bibinfo {title} {Approximation of large-scale dynamical systems}}}\ (\bibinfo  {publisher} {SIAM},\ \bibinfo {year} {2005})\BibitemShut {NoStop}%
\bibitem [{\citenamefont {Benner}\ \emph {et~al.}(2017)\citenamefont {Benner}, \citenamefont {Ohlberger}, \citenamefont {Cohen},\ and\ \citenamefont {Willcox}}]{bennerbook}%
  \BibitemOpen
  \bibfield  {author} {\bibinfo {author} {\bibfnamefont {P.}~\bibnamefont {Benner}}, \bibinfo {author} {\bibfnamefont {M.}~\bibnamefont {Ohlberger}}, \bibinfo {author} {\bibfnamefont {A.}~\bibnamefont {Cohen}},\ and\ \bibinfo {author} {\bibfnamefont {K.}~\bibnamefont {Willcox}},\ }\href@noop {} {\emph {\bibinfo {title} {Model reduction and approximation: theory and algorithms}}}\ (\bibinfo  {publisher} {SIAM},\ \bibinfo {year} {2017})\BibitemShut {NoStop}%
\bibitem [{\citenamefont {LeGresley}\ and\ \citenamefont {Alonso}(2000)}]{legresley2000airfoil}%
  \BibitemOpen
  \bibfield  {author} {\bibinfo {author} {\bibfnamefont {P.}~\bibnamefont {LeGresley}}\ and\ \bibinfo {author} {\bibfnamefont {J.}~\bibnamefont {Alonso}},\ }\bibfield  {title} {\bibinfo {title} {Airfoil design optimization using reduced order models based on proper orthogonal decomposition, 2000},\ }\href@noop {} {\bibfield  {journal} {\bibinfo  {journal} {American Institute of Aeronautics and Astronautics Paper}\ }\textbf {\bibinfo {volume} {2545}} (\bibinfo {year} {2000})}\BibitemShut {NoStop}%
\bibitem [{\citenamefont {Weickum}\ \emph {et~al.}(2009)\citenamefont {Weickum}, \citenamefont {Eldred},\ and\ \citenamefont {Maute}}]{weickum2009multi}%
  \BibitemOpen
  \bibfield  {author} {\bibinfo {author} {\bibfnamefont {G.}~\bibnamefont {Weickum}}, \bibinfo {author} {\bibfnamefont {M.}~\bibnamefont {Eldred}},\ and\ \bibinfo {author} {\bibfnamefont {K.}~\bibnamefont {Maute}},\ }\bibfield  {title} {\bibinfo {title} {A multi-point reduced-order modeling approach of transient structural dynamics with application to robust design optimization},\ }\href@noop {} {\bibfield  {journal} {\bibinfo  {journal} {Structural and Multidisciplinary Optimization}\ }\textbf {\bibinfo {volume} {38}},\ \bibinfo {pages} {599} (\bibinfo {year} {2009})}\BibitemShut {NoStop}%
\bibitem [{\citenamefont {Amsallem}\ \emph {et~al.}(2015)\citenamefont {Amsallem}, \citenamefont {Zahr}, \citenamefont {Choi},\ and\ \citenamefont {Farhat}}]{amsallem2015design}%
  \BibitemOpen
  \bibfield  {author} {\bibinfo {author} {\bibfnamefont {D.}~\bibnamefont {Amsallem}}, \bibinfo {author} {\bibfnamefont {M.}~\bibnamefont {Zahr}}, \bibinfo {author} {\bibfnamefont {Y.}~\bibnamefont {Choi}},\ and\ \bibinfo {author} {\bibfnamefont {C.}~\bibnamefont {Farhat}},\ }\bibfield  {title} {\bibinfo {title} {Design optimization using hyper-reduced-order models},\ }\href@noop {} {\bibfield  {journal} {\bibinfo  {journal} {Structural and Multidisciplinary Optimization}\ }\textbf {\bibinfo {volume} {51}},\ \bibinfo {pages} {919} (\bibinfo {year} {2015})}\BibitemShut {NoStop}%
\bibitem [{\citenamefont {Amsallem}\ and\ \citenamefont {Farhat}(2008)}]{amsallem2008aiaa}%
  \BibitemOpen
  \bibfield  {author} {\bibinfo {author} {\bibfnamefont {D.}~\bibnamefont {Amsallem}}\ and\ \bibinfo {author} {\bibfnamefont {C.}~\bibnamefont {Farhat}},\ }\bibfield  {title} {\bibinfo {title} {Interpolation method for adapting reduced-order models and application to aeroelasticity},\ }\href@noop {} {\bibfield  {journal} {\bibinfo  {journal} {AIAA Journal}\ }\textbf {\bibinfo {volume} {46}},\ \bibinfo {pages} {1803} (\bibinfo {year} {2008})}\BibitemShut {NoStop}%
\bibitem [{\citenamefont {Degroote}\ \emph {et~al.}(2010)\citenamefont {Degroote}, \citenamefont {Vierendeels},\ and\ \citenamefont {Willcox}}]{degroote2010interpolation}%
  \BibitemOpen
  \bibfield  {author} {\bibinfo {author} {\bibfnamefont {J.}~\bibnamefont {Degroote}}, \bibinfo {author} {\bibfnamefont {J.}~\bibnamefont {Vierendeels}},\ and\ \bibinfo {author} {\bibfnamefont {K.}~\bibnamefont {Willcox}},\ }\bibfield  {title} {\bibinfo {title} {Interpolation among reduced-order matrices to obtain parameterized models for design, optimization and probabilistic analysis},\ }\href@noop {} {\bibfield  {journal} {\bibinfo  {journal} {International Journal for Numerical Methods in Fluids}\ }\textbf {\bibinfo {volume} {63}},\ \bibinfo {pages} {207} (\bibinfo {year} {2010})}\BibitemShut {NoStop}%
\bibitem [{\citenamefont {Amsallem}\ and\ \citenamefont {Farhat}(2011)}]{amsallem2011online}%
  \BibitemOpen
  \bibfield  {author} {\bibinfo {author} {\bibfnamefont {D.}~\bibnamefont {Amsallem}}\ and\ \bibinfo {author} {\bibfnamefont {C.}~\bibnamefont {Farhat}},\ }\bibfield  {title} {\bibinfo {title} {An online method for interpolating linear parametric reduced-order models},\ }\href@noop {} {\bibfield  {journal} {\bibinfo  {journal} {SIAM Journal on Scientific Computing}\ }\textbf {\bibinfo {volume} {33}},\ \bibinfo {pages} {2169} (\bibinfo {year} {2011})}\BibitemShut {NoStop}%
\bibitem [{\citenamefont {Amsallem}\ \emph {et~al.}(2016)\citenamefont {Amsallem}, \citenamefont {Tezaur},\ and\ \citenamefont {Farhat}}]{amsallem2016real}%
  \BibitemOpen
  \bibfield  {author} {\bibinfo {author} {\bibfnamefont {D.}~\bibnamefont {Amsallem}}, \bibinfo {author} {\bibfnamefont {R.}~\bibnamefont {Tezaur}},\ and\ \bibinfo {author} {\bibfnamefont {C.}~\bibnamefont {Farhat}},\ }\bibfield  {title} {\bibinfo {title} {Real-time solution of linear computational problems using databases of parametric reduced-order models with arbitrary underlying meshes},\ }\href@noop {} {\bibfield  {journal} {\bibinfo  {journal} {Journal of Computational Physics}\ }\textbf {\bibinfo {volume} {326}},\ \bibinfo {pages} {373} (\bibinfo {year} {2016})}\BibitemShut {NoStop}%
\bibitem [{\citenamefont {Choi}\ \emph {et~al.}(2020)\citenamefont {Choi}, \citenamefont {Boncoraglio}, \citenamefont {Anderson}, \citenamefont {Amsallem},\ and\ \citenamefont {Farhat}}]{choi2020gradient}%
  \BibitemOpen
  \bibfield  {author} {\bibinfo {author} {\bibfnamefont {Y.}~\bibnamefont {Choi}}, \bibinfo {author} {\bibfnamefont {G.}~\bibnamefont {Boncoraglio}}, \bibinfo {author} {\bibfnamefont {S.}~\bibnamefont {Anderson}}, \bibinfo {author} {\bibfnamefont {D.}~\bibnamefont {Amsallem}},\ and\ \bibinfo {author} {\bibfnamefont {C.}~\bibnamefont {Farhat}},\ }\bibfield  {title} {\bibinfo {title} {Gradient-based constrained optimization using a database of linear reduced-order models},\ }\href@noop {} {\bibfield  {journal} {\bibinfo  {journal} {Journal of Computational Physics}\ }\textbf {\bibinfo {volume} {423}},\ \bibinfo {pages} {109787} (\bibinfo {year} {2020})}\BibitemShut {NoStop}%
\bibitem [{\citenamefont {Hess}\ \emph {et~al.}(2023)\citenamefont {Hess}, \citenamefont {Quaini},\ and\ \citenamefont {Rozza}}]{hess2023data}%
  \BibitemOpen
  \bibfield  {author} {\bibinfo {author} {\bibfnamefont {M.~W.}\ \bibnamefont {Hess}}, \bibinfo {author} {\bibfnamefont {A.}~\bibnamefont {Quaini}},\ and\ \bibinfo {author} {\bibfnamefont {G.}~\bibnamefont {Rozza}},\ }\bibfield  {title} {\bibinfo {title} {A data-driven surrogate modeling approach for time-dependent incompressible navier-stokes equations with dynamic mode decomposition and manifold interpolation},\ }\href@noop {} {\bibfield  {journal} {\bibinfo  {journal} {Advances in Computational Mathematics}\ }\textbf {\bibinfo {volume} {49}},\ \bibinfo {pages} {22} (\bibinfo {year} {2023})}\BibitemShut {NoStop}%
\bibitem [{\citenamefont {T{\'o}th}(2010)}]{toth2010modeling}%
  \BibitemOpen
  \bibfield  {author} {\bibinfo {author} {\bibfnamefont {R.}~\bibnamefont {T{\'o}th}},\ }\href@noop {} {\emph {\bibinfo {title} {Modeling and identification of linear parameter-varying systems}}},\ Vol.\ \bibinfo {volume} {403}\ (\bibinfo  {publisher} {Springer},\ \bibinfo {year} {2010})\BibitemShut {NoStop}%
\bibitem [{\citenamefont {Iannelli}\ \emph {et~al.}(2021)\citenamefont {Iannelli}, \citenamefont {Fasel},\ and\ \citenamefont {Smith}}]{iannelli2021ast}%
  \BibitemOpen
  \bibfield  {author} {\bibinfo {author} {\bibfnamefont {A.}~\bibnamefont {Iannelli}}, \bibinfo {author} {\bibfnamefont {U.}~\bibnamefont {Fasel}},\ and\ \bibinfo {author} {\bibfnamefont {R.~S.}\ \bibnamefont {Smith}},\ }\bibfield  {title} {\bibinfo {title} {The balanced mode decomposition algorithm for data-driven lpv low-order models of aeroservoelastic systems},\ }\href@noop {} {\bibfield  {journal} {\bibinfo  {journal} {Aerospace Science and Technology}\ }\textbf {\bibinfo {volume} {115}},\ \bibinfo {pages} {106821} (\bibinfo {year} {2021})}\BibitemShut {NoStop}%
\bibitem [{\citenamefont {Fonzi}\ \emph {et~al.}(2020)\citenamefont {Fonzi}, \citenamefont {Brunton},\ and\ \citenamefont {Fasel}}]{fonzi2020prsa}%
  \BibitemOpen
  \bibfield  {author} {\bibinfo {author} {\bibfnamefont {N.}~\bibnamefont {Fonzi}}, \bibinfo {author} {\bibfnamefont {S.~L.}\ \bibnamefont {Brunton}},\ and\ \bibinfo {author} {\bibfnamefont {U.}~\bibnamefont {Fasel}},\ }\bibfield  {title} {\bibinfo {title} {Data-driven nonlinear aeroelastic models of morphing wings for control},\ }\href@noop {} {\bibfield  {journal} {\bibinfo  {journal} {Proceedings of the Royal Society A}\ }\textbf {\bibinfo {volume} {476}},\ \bibinfo {pages} {20200079} (\bibinfo {year} {2020})}\BibitemShut {NoStop}%
\bibitem [{\citenamefont {Heinze}\ and\ \citenamefont {Borglund}(2008)}]{heinze2008robust}%
  \BibitemOpen
  \bibfield  {author} {\bibinfo {author} {\bibfnamefont {S.}~\bibnamefont {Heinze}}\ and\ \bibinfo {author} {\bibfnamefont {D.}~\bibnamefont {Borglund}},\ }\bibfield  {title} {\bibinfo {title} {Robust flutter analysis considering mode shape variations},\ }\href@noop {} {\bibfield  {journal} {\bibinfo  {journal} {Journal of Aircraft}\ }\textbf {\bibinfo {volume} {45}},\ \bibinfo {pages} {1070} (\bibinfo {year} {2008})}\BibitemShut {NoStop}%
\bibitem [{\citenamefont {Sorokin}\ and\ \citenamefont {Thomsen}(2015)}]{sorokin2015eigenfrequencies}%
  \BibitemOpen
  \bibfield  {author} {\bibinfo {author} {\bibfnamefont {V.~S.}\ \bibnamefont {Sorokin}}\ and\ \bibinfo {author} {\bibfnamefont {J.~J.}\ \bibnamefont {Thomsen}},\ }\bibfield  {title} {\bibinfo {title} {Eigenfrequencies and eigenmodes of a beam with periodically continuously varying spatial properties},\ }\href@noop {} {\bibfield  {journal} {\bibinfo  {journal} {Journal of Sound and Vibration}\ }\textbf {\bibinfo {volume} {347}},\ \bibinfo {pages} {14} (\bibinfo {year} {2015})}\BibitemShut {NoStop}%
\bibitem [{\citenamefont {Xing}\ and\ \citenamefont {Wang}(2024)}]{xing2024dynamic}%
  \BibitemOpen
  \bibfield  {author} {\bibinfo {author} {\bibfnamefont {W.~C.}\ \bibnamefont {Xing}}\ and\ \bibinfo {author} {\bibfnamefont {Y.~Q.}\ \bibnamefont {Wang}},\ }\bibfield  {title} {\bibinfo {title} {Dynamic modeling and vibration analysis of bolted flange joint disk-drum structures: Theory and experiment},\ }\href@noop {} {\bibfield  {journal} {\bibinfo  {journal} {International Journal of Mechanical Sciences}\ }\textbf {\bibinfo {volume} {272}},\ \bibinfo {pages} {109186} (\bibinfo {year} {2024})}\BibitemShut {NoStop}%
\bibitem [{\citenamefont {Gallina}\ \emph {et~al.}(2011)\citenamefont {Gallina}, \citenamefont {Pichler},\ and\ \citenamefont {Uhl}}]{gallina2011enhanced}%
  \BibitemOpen
  \bibfield  {author} {\bibinfo {author} {\bibfnamefont {A.}~\bibnamefont {Gallina}}, \bibinfo {author} {\bibfnamefont {L.}~\bibnamefont {Pichler}},\ and\ \bibinfo {author} {\bibfnamefont {T.}~\bibnamefont {Uhl}},\ }\bibfield  {title} {\bibinfo {title} {Enhanced meta-modelling technique for analysis of mode crossing, mode veering and mode coalescence in structural dynamics},\ }\href@noop {} {\bibfield  {journal} {\bibinfo  {journal} {Mechanical Systems and Signal Processing}\ }\textbf {\bibinfo {volume} {25}},\ \bibinfo {pages} {2297} (\bibinfo {year} {2011})}\BibitemShut {NoStop}%
\bibitem [{\citenamefont {Lu}\ \emph {et~al.}(2019)\citenamefont {Lu}, \citenamefont {Zhan}, \citenamefont {Apley},\ and\ \citenamefont {Chen}}]{lu2019uncertainty}%
  \BibitemOpen
  \bibfield  {author} {\bibinfo {author} {\bibfnamefont {J.}~\bibnamefont {Lu}}, \bibinfo {author} {\bibfnamefont {Z.}~\bibnamefont {Zhan}}, \bibinfo {author} {\bibfnamefont {D.~W.}\ \bibnamefont {Apley}},\ and\ \bibinfo {author} {\bibfnamefont {W.}~\bibnamefont {Chen}},\ }\bibfield  {title} {\bibinfo {title} {Uncertainty propagation of frequency response functions using a multi-output gaussian process model},\ }\href@noop {} {\bibfield  {journal} {\bibinfo  {journal} {Computers \& Structures}\ }\textbf {\bibinfo {volume} {217}},\ \bibinfo {pages} {1} (\bibinfo {year} {2019})}\BibitemShut {NoStop}%
\bibitem [{\citenamefont {Gibanica}\ and\ \citenamefont {Abrahamsson}(2021)}]{gibanica2021data}%
  \BibitemOpen
  \bibfield  {author} {\bibinfo {author} {\bibfnamefont {M.}~\bibnamefont {Gibanica}}\ and\ \bibinfo {author} {\bibfnamefont {T.~J.}\ \bibnamefont {Abrahamsson}},\ }\bibfield  {title} {\bibinfo {title} {Data-driven modal surrogate model for frequency response uncertainty propagation},\ }\href@noop {} {\bibfield  {journal} {\bibinfo  {journal} {Probabilistic Engineering Mechanics}\ }\textbf {\bibinfo {volume} {66}},\ \bibinfo {pages} {103142} (\bibinfo {year} {2021})}\BibitemShut {NoStop}%
\bibitem [{\citenamefont {Lee}\ \emph {et~al.}(2024)\citenamefont {Lee}, \citenamefont {Park}, \citenamefont {Lim},\ and\ \citenamefont {Park}}]{lee2024model}%
  \BibitemOpen
  \bibfield  {author} {\bibinfo {author} {\bibfnamefont {G.-Y.}\ \bibnamefont {Lee}}, \bibinfo {author} {\bibfnamefont {K.-J.}\ \bibnamefont {Park}}, \bibinfo {author} {\bibfnamefont {D.-G.}\ \bibnamefont {Lim}},\ and\ \bibinfo {author} {\bibfnamefont {Y.-H.}\ \bibnamefont {Park}},\ }\bibfield  {title} {\bibinfo {title} {Model order reduction based on low-rank approximation for parameterized eigenvalue problems in structural dynamics},\ }\href@noop {} {\bibfield  {journal} {\bibinfo  {journal} {Journal of Sound and Vibration}\ }\textbf {\bibinfo {volume} {582}},\ \bibinfo {pages} {118413} (\bibinfo {year} {2024})}\BibitemShut {NoStop}%
\bibitem [{\citenamefont {Quarteroni}\ \emph {et~al.}(2015)\citenamefont {Quarteroni}, \citenamefont {Manzoni},\ and\ \citenamefont {Negri}}]{rbbook}%
  \BibitemOpen
  \bibfield  {author} {\bibinfo {author} {\bibfnamefont {A.}~\bibnamefont {Quarteroni}}, \bibinfo {author} {\bibfnamefont {A.}~\bibnamefont {Manzoni}},\ and\ \bibinfo {author} {\bibfnamefont {F.}~\bibnamefont {Negri}},\ }\href@noop {} {\emph {\bibinfo {title} {Reduced basis methods for partial differential equations: an introduction}}}\ (\bibinfo  {publisher} {Springer},\ \bibinfo {year} {2015})\BibitemShut {NoStop}%
\bibitem [{\citenamefont {Antoulas}\ \emph {et~al.}(2020)\citenamefont {Antoulas}, \citenamefont {Beattie},\ and\ \citenamefont {G{\"u}{\u{g}}ercin}}]{loewnerbook}%
  \BibitemOpen
  \bibfield  {author} {\bibinfo {author} {\bibfnamefont {A.~C.}\ \bibnamefont {Antoulas}}, \bibinfo {author} {\bibfnamefont {C.~A.}\ \bibnamefont {Beattie}},\ and\ \bibinfo {author} {\bibfnamefont {S.}~\bibnamefont {G{\"u}{\u{g}}ercin}},\ }\href@noop {} {\emph {\bibinfo {title} {Interpolatory methods for model reduction}}}\ (\bibinfo  {publisher} {SIAM},\ \bibinfo {year} {2020})\BibitemShut {NoStop}%
\bibitem [{\citenamefont {Brunton}\ and\ \citenamefont {Kutz}(2019)}]{databook}%
  \BibitemOpen
  \bibfield  {author} {\bibinfo {author} {\bibfnamefont {S.~L.}\ \bibnamefont {Brunton}}\ and\ \bibinfo {author} {\bibfnamefont {J.~N.}\ \bibnamefont {Kutz}},\ }\href@noop {} {\emph {\bibinfo {title} {Data-driven science and engineering: Machine learning, dynamical systems, and control}}}\ (\bibinfo  {publisher} {Cambridge University Press},\ \bibinfo {year} {2019})\BibitemShut {NoStop}%
\bibitem [{\citenamefont {Gavish}\ and\ \citenamefont {Donoho}(2014)}]{gavish2014optimal}%
  \BibitemOpen
  \bibfield  {author} {\bibinfo {author} {\bibfnamefont {M.}~\bibnamefont {Gavish}}\ and\ \bibinfo {author} {\bibfnamefont {D.~L.}\ \bibnamefont {Donoho}},\ }\bibfield  {title} {\bibinfo {title} {The optimal hard threshold for singular values is $4/\sqrt{3}$},\ }\href@noop {} {\bibfield  {journal} {\bibinfo  {journal} {IEEE Transactions on Information Theory}\ }\textbf {\bibinfo {volume} {60}},\ \bibinfo {pages} {5040} (\bibinfo {year} {2014})}\BibitemShut {NoStop}%
\bibitem [{\citenamefont {Berkooz}\ \emph {et~al.}(1993)\citenamefont {Berkooz}, \citenamefont {Holmes},\ and\ \citenamefont {Lumley}}]{berkooz1993arfm}%
  \BibitemOpen
  \bibfield  {author} {\bibinfo {author} {\bibfnamefont {G.}~\bibnamefont {Berkooz}}, \bibinfo {author} {\bibfnamefont {P.}~\bibnamefont {Holmes}},\ and\ \bibinfo {author} {\bibfnamefont {J.~L.}\ \bibnamefont {Lumley}},\ }\bibfield  {title} {\bibinfo {title} {The proper orthogonal decomposition in the analysis of turbulent flows},\ }\href@noop {} {\bibfield  {journal} {\bibinfo  {journal} {Annual Review of Fluid Mechanics}\ }\textbf {\bibinfo {volume} {25}},\ \bibinfo {pages} {539} (\bibinfo {year} {1993})}\BibitemShut {NoStop}%
\bibitem [{\citenamefont {Gupta}\ \emph {et~al.}(2023)\citenamefont {Gupta}, \citenamefont {Zhao}, \citenamefont {Sharma}, \citenamefont {Agrawal}, \citenamefont {Hourigan},\ and\ \citenamefont {Thompson}}]{gupta2023jfm}%
  \BibitemOpen
  \bibfield  {author} {\bibinfo {author} {\bibfnamefont {S.}~\bibnamefont {Gupta}}, \bibinfo {author} {\bibfnamefont {J.}~\bibnamefont {Zhao}}, \bibinfo {author} {\bibfnamefont {A.}~\bibnamefont {Sharma}}, \bibinfo {author} {\bibfnamefont {A.}~\bibnamefont {Agrawal}}, \bibinfo {author} {\bibfnamefont {K.}~\bibnamefont {Hourigan}},\ and\ \bibinfo {author} {\bibfnamefont {M.~C.}\ \bibnamefont {Thompson}},\ }\bibfield  {title} {\bibinfo {title} {Two-and three-dimensional wake transitions of a naca0012 airfoil},\ }\href@noop {} {\bibfield  {journal} {\bibinfo  {journal} {Journal of Fluid Mechanics}\ }\textbf {\bibinfo {volume} {954}},\ \bibinfo {pages} {A26} (\bibinfo {year} {2023})}\BibitemShut {NoStop}%
\bibitem [{\citenamefont {Fischer}\ \emph {et~al.}(2008)\citenamefont {Fischer}, \citenamefont {Lottes},\ and\ \citenamefont {Kerkemeier}}]{fischer2008nek5000}%
  \BibitemOpen
  \bibfield  {author} {\bibinfo {author} {\bibfnamefont {P.}~\bibnamefont {Fischer}}, \bibinfo {author} {\bibfnamefont {J.}~\bibnamefont {Lottes}},\ and\ \bibinfo {author} {\bibfnamefont {S.}~\bibnamefont {Kerkemeier}},\ }\bibfield  {title} {\bibinfo {title} {Nek5000: open source spectral element cfd solver (2008)},\ }\href@noop {} {\bibfield  {journal} {\bibinfo  {journal} {URL http://nek5000. mcs. anl. gov/index. php/MainPage}\ } (\bibinfo {year} {2008})}\BibitemShut {NoStop}%
\bibitem [{\citenamefont {Frantz}\ \emph {et~al.}(2023)\citenamefont {Frantz}, \citenamefont {Loiseau},\ and\ \citenamefont {Robinet}}]{frantz2023amr}%
  \BibitemOpen
  \bibfield  {author} {\bibinfo {author} {\bibfnamefont {R.~A.~S.}\ \bibnamefont {Frantz}}, \bibinfo {author} {\bibfnamefont {J.-C.}\ \bibnamefont {Loiseau}},\ and\ \bibinfo {author} {\bibfnamefont {J.-C.}\ \bibnamefont {Robinet}},\ }\bibfield  {title} {\bibinfo {title} {{Krylov Methods for Large-Scale Dynamical Systems: Application in Fluid Dynamics}},\ }\bibfield  {journal} {\bibinfo  {journal} {Applied Mechanics Reviews}\ }\textbf {\bibinfo {volume} {75}},\ \href {https://doi.org/10.1115/1.4056808} {10.1115/1.4056808} (\bibinfo {year} {2023}),\ \bibinfo {note} {030802}\BibitemShut {NoStop}%
\bibitem [{\citenamefont {{\AA}kervik}\ \emph {et~al.}(2006)\citenamefont {{\AA}kervik}, \citenamefont {Brandt}, \citenamefont {Henningson}, \citenamefont {H{\oe}pffner}, \citenamefont {Marxen},\ and\ \citenamefont {Schlatter}}]{aakervik2006pof}%
  \BibitemOpen
  \bibfield  {author} {\bibinfo {author} {\bibfnamefont {E.}~\bibnamefont {{\AA}kervik}}, \bibinfo {author} {\bibfnamefont {L.}~\bibnamefont {Brandt}}, \bibinfo {author} {\bibfnamefont {D.~S.}\ \bibnamefont {Henningson}}, \bibinfo {author} {\bibfnamefont {J.}~\bibnamefont {H{\oe}pffner}}, \bibinfo {author} {\bibfnamefont {O.}~\bibnamefont {Marxen}},\ and\ \bibinfo {author} {\bibfnamefont {P.}~\bibnamefont {Schlatter}},\ }\bibfield  {title} {\bibinfo {title} {Steady solutions of the navier-stokes equations by selective frequency damping},\ }\href@noop {} {\bibfield  {journal} {\bibinfo  {journal} {Physics of fluids}\ }\textbf {\bibinfo {volume} {18}} (\bibinfo {year} {2006})}\BibitemShut {NoStop}%
\bibitem [{\citenamefont {Lu}\ \emph {et~al.}(2020)\citenamefont {Lu}, \citenamefont {Tang}, \citenamefont {Apley}, \citenamefont {Zhan},\ and\ \citenamefont {Chen}}]{lu2020mode}%
  \BibitemOpen
  \bibfield  {author} {\bibinfo {author} {\bibfnamefont {J.}~\bibnamefont {Lu}}, \bibinfo {author} {\bibfnamefont {J.}~\bibnamefont {Tang}}, \bibinfo {author} {\bibfnamefont {D.~W.}\ \bibnamefont {Apley}}, \bibinfo {author} {\bibfnamefont {Z.}~\bibnamefont {Zhan}},\ and\ \bibinfo {author} {\bibfnamefont {W.}~\bibnamefont {Chen}},\ }\bibfield  {title} {\bibinfo {title} {A mode tracking method in modal metamodeling for structures with clustered eigenvalues},\ }\href@noop {} {\bibfield  {journal} {\bibinfo  {journal} {Computer Methods in Applied Mechanics and Engineering}\ }\textbf {\bibinfo {volume} {369}},\ \bibinfo {pages} {113174} (\bibinfo {year} {2020})}\BibitemShut {NoStop}%
\end{thebibliography}
\end{document}